\documentclass[11pt,a4paper,reqno]{amsart}
\usepackage{amssymb,amsmath,mathrsfs,enumerate}
\usepackage[margin=3cm]{geometry}
\usepackage{color}
\usepackage{hyperref}
\hypersetup{colorlinks=true,urlcolor=blue,linkcolor=blue,citecolor=red,bookmarksopen=true}

\numberwithin{equation}{section}
\theoremstyle{plain}
\newtheorem{Prop}{Proposition}[section]
\newtheorem{Thm}[Prop]{Theorem}
\newtheorem{Lem}[Prop]{Lemma}
\newtheorem{Cor}[Prop]{Corollary}
\theoremstyle{definition}
\newtheorem{Def}[Prop]{Definition}

\newtheorem{Prob}[Prop]{Open problem}
\theoremstyle{remark}
\newtheorem{Rem}[Prop]{Remark}
 \newtheorem*{ThmA}{Theorem A}

\DeclareMathOperator{\spn}{span}
\DeclareMathOperator{\diam}{diam}
\DeclareMathOperator{\dv}{div}
\DeclareMathOperator{\Tr}{Tr}
\DeclareMathOperator{\dist}{dist}

\DeclareMathOperator{\argmax}{argmax}
\DeclareMathOperator{\argmin}{argmin}
\DeclareMathOperator{\USC}{USC}
\DeclareMathOperator{\LSC}{LSC}
\DeclareMathOperator{\osc}{osc}

\newcommand{\emp}{\emptyset}

\newcommand{\mns}{\setminus}
\newcommand{\N}{\mathbb{N}}
\newcommand{\R}{\mathbb{R}}
\newcommand{\G}{\mathbb{G}}
\newcommand{\fr}[1]{\mathfrak{#1}}
\newcommand{\g}{\fr{g}}
\newcommand{\del}{\partial}
\newcommand{\I}{\mathbb{I}}
\newcommand{\eps}{\varepsilon}
\newcommand{\inv}[1]{{#1}^{-1}}

\newcommand{\Om}{\Omega}
\newcommand{\om}{\omega}
\newcommand{\inp}[2]{\big\langle #1,#2\big\rangle}
\newcommand{\gr}{\nabla}
\newcommand{\lap}{\Delta}

\newcommand{\hh}{\mathbb{H}}
\newcommand{\X}{\mathfrak{X}}

\newcommand{\XX}{\X\X}
\newcommand{\XXs}{\XX^\star}
\newcommand{\dvh}{\dv_{{\mathstrut}\X}}
\newcommand{\A}{\mathcal{A}}
\newcommand{\Av}{\mathscr{A}}

\newcommand{\E}{\mathcal{E}}
\newcommand{\LL}{\mathcal{L}}
\newcommand{\cc}{\mathcal{C}}
\newcommand{\Z}{\mathcal{Z}}
\newcommand{\hxi}{\hat\xi}

\title[Comparison principles for sub-elliptic equations]
{Comparison principles for degenerate sub-elliptic equations in non-divergence form}

\author[Juan J. Manfredi]{Juan J. Manfredi}
\address[J.J.\ Manfredi]{Department of Mathematics,
University of Pittsburgh, 312 Thackeray Hall, Pittsburgh, PA 15260, USA.}
\email{manfredi@pitt.edu}

\author[Shirsho Mukherjee]{Shirsho Mukherjee}
\address[S.\ Mukherjee]{Department of Mathematical Sciences,
University of Essex, Wivenhoe Park, Colchester, CO4 3SQ, UK.}
\email{shirsho.mukherjee@essex.ac.uk}
\date{October 18, 2026}
\begin{document}

 \begin{abstract}
We study the comparison principle for viscosity sub- and super-solutions of degenerate quasi-linear equations in non-divergence form $-\Tr\big(A(\X u)\,\XXs u\big)+H(\X u)=0$, driven by the horizontal vector fields of a Carnot group, under the ellipticity condition of Barles and Busca and a homogeneity condition along rays. The class contains the sub-elliptic $\infty$-Laplacian and the normalized $p$-Laplacian, and also operators of rank one that are not Aronsson operators, such as the tilted $\infty$-Laplacian $\inp{\XXs u\,T\X u}{T\X u}$ with $T$ not symmetric. We use the strong maximum principle of Bardi and Goffi to implement the Barles and Busca method, and we identify the obstruction created by non-commutativity: the translations that preserve the equation and the translations compatible with the first-order Taylor expansion in the horizontal gradient are different, and the only ones that do both are the central ones. We prove the comparison principle for pairs of functions whose regularizations by sup- and inf-convolution have non-vanishing horizontal gradient on the maximal set of their difference, and unconditionally when one of the two functions is $C^1$ with non-vanishing horizontal gradient near that set, when one of them is constant, and when both depend only on the horizontal coordinates, which includes the Euclidean case. Whether the hypothesis can be removed is stated as an open problem.
\end{abstract}

\maketitle

\section{Introduction}\label{sec:Introduction}

Consider the Dirichlet problem for the $\infty$-Laplacian in a bounded domain $\Omega\subset\mathbb{R}^n$ with $f\in C^{\,0,1}(\del\Om)$,
\begin{equation}\label{eq:infinitylap0}
\begin{cases}
\lap_\infty u  =  0\quad &\text{in}\ \ \Om\\
 \quad \ \ u  =  f \quad &\text{on}\ \, \del\Om,
\end{cases}
\end{equation}
where the $\infty$-Laplacian $\lap_\infty u =\inp{D^2 u\gr u}{\gr u}$ is a degenerate elliptic operator and the solutions of \eqref{eq:infinitylap0} are considered in the viscosity sense. The uniqueness of viscosity solutions of \eqref{eq:infinitylap0} was established in a celebrated theorem of Jensen \cite{J93}, based on the approximation of the $\infty$-Laplacian by the $p$-Laplacian as $p\to\infty$ and on auxiliary upper and lower envelope equations.

Uniqueness is a consequence of the \emph{comparison principle}: if $u$ is a viscosity sub-solution and $v$ is a viscosity super-solution in $\Om$, and $u\leq v$ on $\del\Om$, then $u\leq v$ in $\Om$. In the standard comparison theorems for viscosity solutions \cite{Cr-Ish-Lions} the equation is assumed to be strictly monotone in the unknown $u$, and this assumption is not available for equations in which, as in $\lap_\infty u=0$, the unknown appears only through its derivatives. Barles and Busca \cite{Barles-Busca} gave a proof of the comparison principle for a general class of degenerate elliptic equations of this kind, which includes $\Delta_\infty u=0$ and equations of the form
\begin{equation}\label{eq:bb}
-\Tr \big(A(\gr u)D^2u\big)+ H(\gr u)=0\qquad \text{in}\ \Om,
\end{equation}
where $A(\xi)$ is a symmetric and positive semidefinite matrix. Their approach does not use approximations by $p$-Laplacians and is based on the strong maximum principle; the equation is not required to have a variational structure. Their conditions on $A$ and $H$ are stated in Section~\ref{subsec:BB}. \par
Using the characterisation of $\infty$-harmonic functions as the functions satisfying comparison with cones, Armstrong and Smart gave a short proof of uniqueness \cite{AS10}; in fact, they prove the uniqueness of absolutely minimizing Lipschitz extensions with given boundary values. Their proof is valid in those settings where $\infty$-harmonic functions are absolute minimizers, or equivalently satisfy comparison with cones \cite{CdP07}. \par
For Riemannian and sub-Riemannian manifolds, absolute minimizers are $\infty$-harmonic functions in the viscosity sense. The converse holds for Riemannian manifolds and for those sub-Riemannian manifolds that are Carnot groups (see, for example, \cite[Theorem 4.10]{DMV13}). Therefore the Dirichlet problem \eqref{eq:infinitylap0} has a unique viscosity solution in Riemannian manifolds and in Carnot groups. For general sub-Riemannian manifolds it is not known whether $\infty$-harmonic functions are absolute minimizers, and the uniqueness for \eqref{eq:infinitylap0} is open.
Bieske \cite{B02} extended Jensen's proof to the Heisenberg group, and C.~Wang \cite{CYW07} extended it to Carnot groups and to the Aronsson equations of supremal functionals; see also \cite{B05, Och-Ruiz}. The approximation method is tied to the variational origin of the equation: the operators it is designed for are the Aronsson operators of convex functions,  whose matrix $A(\xi)=\gr f(\xi)\otimes\gr f(\xi)$ has rank one.
\par
The purpose of this paper is to carry the Barles--Busca method over to sub-elliptic analogues of \eqref{eq:bb} on Carnot groups. The method requires the strong maximum principle for degenerate sub-elliptic equations, recently obtained by Bardi and Goffi \cite{Bardi-Goffi}, Rademacher's theorem on Carnot groups \cite{P89}, and the sup and inf-convolutions of \cite{CYW07}. Every ingredient of the Euclidean proof transfers except one: translations. On a non-abelian group the two roles played by the Euclidean translation $x\mapsto x+h$  are a symmetry of the equation (left translations), and  the first-order Taylor expansion in the gradient (right translations). We prove the comparison principle under a non-degeneracy hypothesis on the regularizations of the pair of functions, which is the hypothesis of the first step of the Barles--Busca argument (Theorem~\ref{thm:comp}), and unconditionally in several cases (Theorem~\ref{thm:uncond}). Whether the hypothesis can be removed is an open problem, which we formulate in Section~\ref{sec:obstruction}.

The rest of the introduction is organized as follows. In Section~\ref{subsec:BB} we state the conditions and the theorem of Barles and Busca in the Euclidean case. In Section~\ref{subsec:defs} we describe the sub-elliptic setting and our hypotheses, and we explain how they are related to the conditions of Barles and Busca. In Section~\ref{subsec:examples} we present the examples, which are collected in Table~\ref{tab:examples}: for each of them the table records the conditions of Barles and Busca, which give uniqueness in the Euclidean case, our hypotheses, which are the ones we use on Carnot groups, and what is known about uniqueness on Carnot groups. The results are stated in Section~\ref{subsec:results}, and Section~\ref{subsec:ideas} describes the obstacle created by non-commutativity and the ideas of the proofs.

\subsection{The conditions of Barles and Busca in the Euclidean case}\label{subsec:BB}
Let $\Om\subset\R^n$ be a bounded domain and consider equation \eqref{eq:bb}, where $A$ is a function on $\R^n$ with values in the symmetric $n\times n$ matrices and $H$ is a real function on $\R^n$. We write $\xi$ for the variable of $A$ and $H$, whose place in the equation is taken by the gradient. The conditions of Barles and Busca are the following three:
\begin{align}
\label{eq:BB1}
\text{(BB1)}&\quad \text{$A$ and $H$ are locally Lipschitz continuous on $\R^n$, and $A(\xi)\geq 0$ for all $\xi$;}\\
\label{eq:BB2}
\text{(BB2)}&\quad \inp{A(\xi)\,\xi}{\xi}>0\quad \text{for all}\ \xi\in \R^n\mns\{0\};\\
\label{eq:BB3}
\text{(BB3)}&\quad \inp{A(\xi)\,\xi}{\xi}\,\geq\, c\,\Big(|\xi|^2\,\Tr A(\xi)+|\xi|\,|H(\xi)|\Big)\quad \text{for}\ 0<|\xi|<\nu,
\end{align}
for some constants $c>0$ and $\nu>0$.
\begin{ThmA}[Barles and Busca {\cite[Theorem 1]{Barles-Busca}}]
Assume \textup{(BB1)}, \textup{(BB2)} and \textup{(BB3)}. Let $u$ be a bounded upper semicontinuous viscosity sub-solution of \eqref{eq:bb} in $\Om$ and let $v$ be a bounded lower semicontinuous viscosity super-solution of \eqref{eq:bb} in $\Om$, such that $u\leq v$ on $\del\Om$. Then $u\leq v$ in $\Om$.
\end{ThmA}
In particular, for given boundary values the Dirichlet problem for \eqref{eq:bb} has at most one viscosity solution in $C(\bar\Om)$.

Condition (BB2) is the key assumption. It is a condition of strict ellipticity in one direction only, the direction of the gradient: the matrix $A(\xi)$ may be degenerate in all the other directions, as it is for the $\infty$-Laplacian, for which $A(\xi)=\xi\otimes\xi$ has rank one. Condition (BB3) concerns small gradients: there, the ellipticity in the direction of the gradient has to be comparable to the size of the whole matrix $A(\xi)$, measured by its trace, and it has to dominate the first-order term. Condition (BB1) expresses the degenerate ellipticity of the equation and the regularity of its coefficients. Note that there is no assumption of growth, of convexity or of homogeneity, and that the equation is not required to come from a variational problem.

The proof of Theorem~A consists of two steps. After a regularization by sup- and inf-convolutions, $u$ and $v$ are differentiable at the points where $u-v$ attains its maximum. In the first step one assumes that the gradient of $u$ does not vanish at these points. A change of unknown $u\mapsto\psi(u)$, $v\mapsto\psi(v)$, with $\psi$ close to the identity, transforms \eqref{eq:bb} into an equation which, by (BB2), is strictly monotone in the unknown as long as the gradient does not vanish; (BB1) is used to control the remaining terms. For such equations the comparison principle follows from  standard arguments. In the second step the assumption on the gradient is removed. Since \eqref{eq:bb} does not depend on $x$, the functions $u(\cdot+h)$ are sub-solutions, and there are two alternatives: either for some small $h$ the gradient does not vanish at the maximum points of $u(\cdot+h)-v$, and the first step applies, or the maximum of $u(\cdot+h)-v$ does not depend on $h$ for small $h$, and then $u$ has a local maximum at every maximum point of $u-v$. Condition (BB3) rules out the second alternative: it implies a Hopf type lemma \cite[Lemma~3.2]{Barles-Busca}, by which a sub-solution is constant near an interior local maximum and a super-solution is constant near an interior local minimum, and this leads to a contradiction with the boundary condition. We come back to this argument in Section~\ref{subsec:ideas}, because the second step is the one that does not extend to Carnot groups.

\emph{The fully non-linear form of the conditions.} The theorem of \cite{Barles-Busca} is stated for equations $F(\gr u,D^2u)=0$, where $F$ is a function of $(\xi,M)\in\R^n\times\mathcal S_n$ and $\mathcal S_n$ is the space of symmetric $n\times n$ matrices, under the following assumptions:
\begin{enumerate}
\item[(F1)] $F$ is locally Lipschitz continuous and degenerate elliptic, that is, $F(\xi,M)\leq F(\xi,N)$ whenever $M\geq N$;
\item[(F2)] for every $\nu>0$ there is $\eps(\nu)>0$ such that $-\Tr\big(F_M(\xi,M)\,\xi\otimes\xi\big)\geq\eps(\nu)$ for almost every $(\xi,M)$ with $\nu\leq|\xi|\leq\nu^{-1}$ and $|M|\leq\nu^{-1}$, where $F_M$ denotes the differential of $F$ with respect to $M$;
\item[(F3)] for every $\alpha>0$ there are $\beta>0$ and $\nu>0$ such that, for $0<|\xi|<\nu$,
\[
F\Big(\xi,\,-\alpha|\xi|\,\I_n+\beta\,\frac{\xi\otimes\xi}{|\xi|}\Big)<0\qquad\text{and}\qquad F\Big(\xi,\,\alpha|\xi|\,\I_n-\beta\,\frac{\xi\otimes\xi}{|\xi|}\Big)>0 .
\]
\end{enumerate}
For $F(\xi,M)=-\Tr\big(A(\xi)M\big)+H(\xi)$ we have $F_M=-A(\xi)$, so that (F1) is (BB1) and, $A$ being continuous, (F2) is (BB2); and the two inequalities in (F3) together read $\beta\inp{A(\xi)\xi}{\xi}>\alpha|\xi|^2\,\Tr A(\xi)+|\xi|\,|H(\xi)|$, so that, given (BB2), (F3) is equivalent to (BB3); see \cite[Section~2]{Barles-Busca}. The examples considered in \cite{Barles-Busca} are the $\infty$-Laplacian, the $p$-Laplacian and the minimal surface equation; we come back to them in Section~\ref{subsec:examples}.

\subsection{The sub-elliptic setting and our hypotheses}\label{subsec:defs}
Let $\G=(\R^n,*)$ be a Carnot group of step $r\geq1$ and let $X_1,\ldots,X_m$ be a basis of the first layer of its Lie algebra of left-invariant vector fields; in particular $X_1,\ldots,X_m$ are $1$-homogeneous with respect to the dilations of $\G$ and satisfy H\"ormander's rank condition
\begin{equation}\label{eq:horcond}
\dim\Big(\texttt{Lie}\big(X_1(x),\ldots,X_m(x)\big) \Big)=n, \qquad \text{for all}\ x\in \R^n.
\end{equation}
Conversely, linearly independent, left-invariant, $1$-homogeneous vector fields on a homogeneous Lie group that satisfy \eqref{eq:horcond} are the first layer of a stratification, so that this is the general setting of left-invariant homogeneous H\"ormander systems (see Section~\ref{subsec:carnot}). The case $m=n$ corresponds to $r=1$ and $\G$ being the Euclidean additive group, with $X_i=\del_{x_i}$. The simplest non-abelian example is the Heisenberg group $\hh^1=\R^3$ with the vector fields
\[
X_1=\del_{x_1}-\tfrac12\,x_2\,\del_{x_3},\qquad X_2=\del_{x_2}+\tfrac12\,x_1\,\del_{x_3},
\]
for which $m=2$, $n=3$ and $[X_1,X_2]=\del_{x_3}$.\par
Given a bounded domain $\Om\subset\R^n$, we consider equations of the form
\begin{equation}\label{eq:eq}
\sum_{i,j=1}^m  A_{i,j}(\X u)X_jX_i u=  H(\X u) \qquad \text{in}\ \Om,
\end{equation}
where $\X u = (X_1u, \ldots,X_m  u)$ is the horizontal gradient and $ \XX u = (X_jX_iu)_{i,j}$ the horizontal second derivative. Here $A_{i,j}, H \in C(\R^m\mns\{0\})$, and the $(m\times m)$ matrix $A(\xi)$ with entries $A_{i,j}(\xi)$ is symmetric and positive semidefinite. Hence \eqref{eq:eq} can be re-written as
\begin{equation}\label{eq:maineq}
\LL u:=-\Tr\big(A(\X u) \XXs u\big) + H(\X u)=0,
\end{equation}
where $\XXs u$ is the symmetrized matrix with entries $\frac{1}{2}(X_iX_ju+X_jX_iu)$. In the Euclidean case $\X u=\gr u$ and $\XXs u=D^2u$, and \eqref{eq:maineq} is \eqref{eq:bb}. Viscosity sub- and super-solutions of \eqref{eq:maineq} are defined with test functions whose horizontal gradient does not vanish at the point of contact (Definition~\ref{def:visc}), so that $A$ and $H$ need not be defined at $\xi=0$.

We assume that $A$ and $H$ satisfy the following conditions of strict ellipticity in the gradient direction and of homogeneity along rays:
\begin{align}
\label{eq:ell}
\text{(E)}&\quad \E(\xi):=\inp{A(\xi)\,\xi}{\xi}>0,\quad \forall\ \xi\in \R^m\mns\{0\};\\
\label{eq:scale}
\text{(S)}&\quad A(t\xi)=\kappa(t,\xi)\,A(\xi)\ \ \text{and}\ \ H(t\xi)= t\,\kappa(t,\xi)\,H(\xi),\quad \forall\ t\geq 1,\ \xi\in\R^m\mns\{0\},
\end{align}
for some function $\kappa:[1,\infty)\times(\R^m\mns\{0\})\to(0,\infty)$. The equation \eqref{eq:maineq} is said to be degenerate if $\lim_{\xi\to 0}A(\xi)=0$.

\begin{Rem}\label{rem:scale}
(a) By \eqref{eq:ell}, $\kappa(t,\xi)=\E(t\xi)/(t^2\E(\xi))$ is continuous, and \eqref{eq:scale} says that along every ray $\{t\xi:t>0\}$ the matrix $A$ and the function $H/|\xi|$ change only by the same positive scalar factor: $A(\xi)=a(\xi)A_0(\xi/|\xi|)$ and $H(\xi)=|\xi|\,a(\xi)H_0(\xi/|\xi|)$, with $a>0$ continuous and $A_0,H_0$ continuous on the unit sphere, $\inp{A_0(\om)\om}{\om}>0$; then $\kappa(t,\xi)=a(t\xi)/a(\xi)$. Since $a>0$ and test functions with vanishing horizontal gradient are not used, the viscosity sub- and super-solutions of \eqref{eq:maineq} are determined by $A_0$ and $H_0$: multiplying $\LL$ by a positive function of $\X u$ does not change them. In particular $A$ and $H$ may be replaced by $A_0(\xi/|\xi|)$ and $|\xi|H_0(\xi/|\xi|)$, for which \eqref{eq:scale} holds with $\kappa\equiv1$; we shall do so from Section~\ref{sec:prelim} on, see \eqref{eq:normal}.

(b) With $\phi(1/t;\xi):=1/(t\kappa(t,\xi))$, condition \eqref{eq:scale} is equivalent to the pair of inequalities
\begin{equation}\label{eq:scaleineq}
-\Tr\big(A(t\xi) X \big) \leq  \frac{ -\Tr\big(A(\xi) X\big)}{t \phi(1/t;\xi)},\qquad H(t\xi)\leq  \frac{H(\xi)}{ \phi(1/t;\xi)},\qquad t\ge1,
\end{equation}
required for all symmetric $X\in\R^{m\times m}$, together with the reversed inequality for $H$. Indeed, for an operator affine in $X$ the first inequality for all $X$ forces $A(t\xi)$ to be a positive multiple of $A(\xi)$ (the map $X\mapsto\Tr\big([A(t\xi)-A(\xi)/(t\phi)]X\big)$ is linear and non-negative, hence zero), so that \eqref{eq:scaleineq} is the scaling hypothesis of Bardi--Goffi \cite[(ii) and Remark 2.8]{Bardi-Goffi} for the sub-solution side, and the reversed inequality for $H$ is their hypothesis (ii$'$) for the super-solution side. Both sides are used below (Corollary \ref{cor:stmx}, Lemma~\ref{lem:case1} and Theorem~\ref{thm:uncond}(ii)), which is why $H$ is required to scale exactly. If $\kappa$ does not depend on $\xi$, then $\kappa(t,\xi)=t^\alpha$ for some $\alpha\in\R$ and \eqref{eq:scale} is exact homogeneity: $A(t\xi)=t^\alpha A(\xi)$, $H(t\xi)=t^{\alpha+1}H(\xi)$.

(c) Condition \eqref{eq:scale} excludes operators whose matrix changes shape along rays. For instance the normalized Uhlenbeck--Ural'tseva operators $A(\xi)=\I_m+\big(|\xi|g'(|\xi|)/g(|\xi|)-1\big)(\xi\otimes\xi)/|\xi|^2$, whose eigenvalues along $\xi$ and orthogonally to $\xi$ are $|\xi|g'(|\xi|)/g(|\xi|)$ and $1$, satisfy \eqref{eq:scale} only when $g$ is a power, and the mean-curvature type operators $\Delta_\X u-\Delta_{\X,\infty}u/(1+|\X u|^2)$, with eigenvalues $1/(1+|\xi|^2)$ and $1$, do not satisfy it.
\end{Rem}

Let us comment on the relationship between the conditions (BB1)--(BB3) and conditions (E) and (S). The latter are conditions on the functions $A$ and $H$ of the gradient variable $\xi$ alone and they do not involve the group. With $\R^m$ in place of $\R^n$ in (BB1)--(BB3), they are related as follows:
\begin{enumerate}[(a)]
\item Condition (E) is (BB2).
\item Condition (S) is not one of the conditions of \cite{Barles-Busca}, where no homogeneity is assumed. It is the form taken, for the operators \eqref{eq:maineq}, by the scaling hypothesis of the strong maximum principle of Bardi and Goffi \cite{Bardi-Goffi}; see Remark~\ref{rem:scale}(b). On a Carnot group this strong maximum principle takes the place of the Hopf lemma that Barles and Busca derive from (BB3).
\item Conditions (E) and (S) imply (BB3), for all $\xi\neq0$. Indeed, by Remark~\ref{rem:scale}(a), with $\hxi=\xi/|\xi|$,
\[
\E(\xi)=a(\xi)\,|\xi|^2\,\E_0(\hxi),\qquad |\xi|^2\,\Tr A(\xi)+|\xi|\,|H(\xi)|=a(\xi)\,|\xi|^2\,\big(\Tr A_0(\hxi)+|H_0(\hxi)|\big),
\]
where $\E_0(\om)=\inp{A_0(\om)\om}{\om}$, so that (BB3) holds with $c=\min\E_0\big/\max\big(\Tr A_0+|H_0|\big)$, the minimum and the maximum being taken on the unit sphere.
\item Instead of the local Lipschitz continuity in (BB1) we assume only that $A$ and $H$ are continuous away from the origin. This is possible because, under (S), the change of unknown in the first step of the proof is given by an exact identity (Lemma~\ref{lem:tech}). Since only the values of $A$ and $H$ on the unit sphere matter (Remark~\ref{rem:scale}(a)), no condition on their growth near the origin is needed either.
\item Conversely, (BB1)--(BB3) do not imply (S). The minimal surface equation, which is one of the examples of \cite{Barles-Busca}, satisfies (BB1)--(BB3) and does not satisfy (S); see row (4) of Table~\ref{tab:examples}.
\end{enumerate}
Thus, in the Euclidean case, the hypotheses of Theorem~A and ours overlap, and neither set contains the other. If $A$ and $H$ satisfy (E) and (S) and the functions $A_0$ and $H_0$ of Remark~\ref{rem:scale}(a) are Lipschitz continuous on the unit sphere, then the functions $|\xi|^2A_0(\hxi)$ and $|\xi|^3H_0(\hxi)$, extended by zero to $\xi=0$, satisfy (BB1)--(BB3) and define an equation with the same viscosity sub- and super-solutions as \eqref{eq:maineq}, so that Theorem~A applies.

On a Carnot group that is not abelian the situation is different, and Theorem~A cannot be applied to \eqref{eq:maineq}. If we write $X_j=\sigma^j(x)\cdot\gr$ and $\sigma(x)=[\sigma^1(x),\ldots,\sigma^m(x)]$, an $n\times m$ matrix, then $\X u=\sigma(x)^T\gr u$, and \eqref{eq:maineq} is an equation $F(x,\gr u,D^2u)=0$ in $\R^n$; see \eqref{eq:sigrad} and \eqref{eq:FtoG}. This equation depends on $x$, and it does not satisfy (BB2) and (BB3) as conditions on the Euclidean gradient $p=\gr u$: its ellipticity in the direction of $p$ is
\[
-\Tr\big(F_M(x,p,M)\,p\otimes p\big)=\E\big(\sigma(x)^Tp\big),
\]
which is defined and positive only when $\sigma(x)^Tp\neq0$. Since $m<n$, at every point $x$ there are vectors $p\neq0$ with $\sigma(x)^Tp=0$; in $\hh^1$ they are the multiples of $(x_2,-x_1,2)$. In other words, the equation is elliptic in the direction of the horizontal gradient, and it does not see the other components of the gradient. This is why the conditions have to be imposed on $A$ and $H$ as functions of the horizontal gradient, as in (E) and (S), and the Hopf lemma has to be replaced by a strong maximum principle adapted to the vector fields.

The strong maximum principle for linear sub-elliptic equations is due to Bony \cite{B69}. It was generalized to non-linear sub-elliptic equations by Capogna--Zhou \cite{CZ18} and Bardi--Goffi \cite{Bardi-Goffi}. The latter also prove a \emph{strong} comparison principle for equations in Hamilton--Jacobi--Bellman form, which is a different statement from the comparison principle considered here.

\subsection{Examples}\label{subsec:examples}
The examples are collected in Table~\ref{tab:examples}. Each equation is given by its matrix $A$ and its first-order term $H$, as in \eqref{eq:maineq}. Since (BB1)--(BB3), (E) and (S) are conditions on $A$ and $H$ only, each example is checked once, with $\xi\in\R^m$. In the Euclidean case $\G=\R^m$ all the examples satisfy (BB1)--(BB3), so that the Dirichlet problem has at most one viscosity solution by Theorem~A. On a Carnot group the conditions that play their role are (E) and (S), which are the hypotheses of Theorems~\ref{thm:comp} and \ref{thm:uncond} below; they hold for all the examples except the minimal surface equation. The last column records what is known about uniqueness for the Dirichlet problem on non-abelian Carnot groups; its entries are explained in the notes to the table.

\begin{table}[t]
\renewcommand{\arraystretch}{1.45}
\footnotesize
\begin{tabular}{@{}l@{\hspace{5pt}}l@{\hspace{5pt}}l@{\hspace{5pt}}l@{\hspace{5pt}}l@{\hspace{5pt}}l@{}}
\hline
 &  &  & In $\R^m$ & \multicolumn{2}{l}{On Carnot groups}\\[-5pt]
Equation & $A(\xi)$ & $H(\xi)$ & (BB1)--(BB3) & (E), (S) & Uniqueness\\
\hline
(1) Laplacian & $\I_m$ & $0$ & yes & yes & yes \cite{B69}\\
(2) $\infty$-Laplacian & $\xi\otimes\xi$ & $0$ & yes & yes & yes \cite{B02, CYW07, B12, DMV13}\\
(3) $p$-Laplacian & $|\xi|^2\I_m{+}(p{-}2)\,\xi\otimes\xi$ & $0$ & yes & yes & yes \cite{B12}\\
(4) Minimal surfaces & $\I_m-\dfrac{\xi\otimes\xi}{1+|\xi|^2}$ & $0$ & yes & (S) fails & not covered\\
(5) Aronsson, $g$ convex & $g^2\,\gr g\otimes\gr g$ & $0$ & yes & yes & yes$^{*}$ \cite{CYW07}\\
(6) Tilted $\infty$-Laplacian & $T\xi\otimes T\xi$ & $0$ & yes & yes & partial\\
(7) Aronsson, $g$ non-convex & $g^2\,\gr g\otimes\gr g$ & $0$ & yes & yes & partial\\
(8) Orthotropic, $p\geq3$ & $\operatorname{diag}(|\xi_i|^{p-2})$ & $0$ & yes & yes & yes$^{**}$\\
(9) $\infty$-Laplacian, biased & $\xi\otimes\xi$ & $-\beta|\xi|^3$ & yes & yes & partial\\
(10) $\infty$-Laplacian, drift & $\xi\otimes\xi$ & $-\langle b,\xi\rangle|\xi|^2$ & yes & yes & partial\\
\hline
\end{tabular}

\smallskip
\parbox{\textwidth}{\footnotesize
$^{*}$ Under the hypotheses of the uniqueness theorem of \cite{CYW07}, which include the strict convexity of $f=\frac12g^2$.\\[2pt]
$^{**}$ For weak solutions: this is a variational problem, the equation being the Euler--Lagrange equation of the strictly convex functional $\int\sum_i|X_iu|^p\,dx$. For viscosity solutions the result is partial, as in rows (6), (7), (9) and (10).\\[2pt]
\emph{Partial}: uniqueness is known only through the conditional Theorems~\ref{thm:comp} and \ref{thm:uncond} of this paper, that is, for pairs of sub- and super-solutions whose regularizations satisfy the non-degeneracy hypothesis of Theorem~\ref{thm:comp} and in the cases listed in Theorem~\ref{thm:uncond}; to our knowledge the general case is open.\\[2pt]
\emph{Not covered}: condition (S) fails, and our results do not apply.}
\medskip
\caption{The examples. The equation is $-\Tr\big(A(\X u)\,\XXs u\big)+H(\X u)=0$. Fourth column: the conditions of Barles and Busca, which give uniqueness for the Dirichlet problem in $\R^m$ (Theorem~A). Fifth column: conditions (E) and (S), the hypotheses of Theorems~\ref{thm:comp} and \ref{thm:uncond}. Last column: uniqueness for the Dirichlet problem on non-abelian Carnot groups.}\label{tab:examples}
\end{table}

In the table, $g\in C^2(\R^m\mns\{0\})$ is a positive function that is positively homogeneous of degree one, $T$ is an $m\times m$ matrix whose symmetric part is positive definite, $\beta\in\R$ and $b\in\R^m$. Only the values of $A$ and $H$ on the unit sphere matter, because multiplying $\LL$ by a positive function of $\X u$ does not change its viscosity sub- and super-solutions (Remark~\ref{rem:scale}(a)). Each equation is written in a form in which $A$ and $H$ are locally Lipschitz continuous on the whole of $\R^m$, as (BB1) requires. Condition (BB2) is (E), and (BB3) follows from (E) and (S), as we saw in Section~\ref{subsec:defs}; in row (4), (BB3) holds because $\E(\xi)=|\xi|^2\,\Tr A(\xi)\big/\big(m+(m-1)|\xi|^2\big)$. Except in rows (1) and (4) we have $A(0)=0$ and $H(0)=0$, so that test functions with vanishing gradient impose no condition, and the viscosity sub- and super-solutions of \cite{Barles-Busca}, which are defined with all test functions, are those of Definition~\ref{def:visc}.

Rows (1)--(5) are standard. Row (1) is the sub-Laplacian $\lap_\X u = \Tr (\XX u)$; the theory of linear equations defined by vector fields satisfying \eqref{eq:horcond} goes back to Chow \cite{Chow} and H\"ormander \cite{Hor}, see \cite{Folland, Folland-Stein, Roth-Stein, Nag-Stein-Wain}, and the uniqueness for its Dirichlet problem follows from the maximum principle of Bony \cite{B69}. Rows (2) and (3) are the sub-elliptic $\infty$-Laplacian and $p$-Laplacian, $1<p<\infty$,
\begin{equation}\label{eq:inflap}
\lap_{\X,\infty} u= \inp{\XXs u \,\X u}{\X u}\qquad\text{and}\qquad \lap_{\X,p} u=\dvh \big(|\X u|^{p-2}\X u\big),
\end{equation}
where $\dvh F= \sum_{i=1}^m X_i F_i$ is the horizontal divergence. We have $\lap_{\X,p} u=|\X u|^{p-2}\lap^N_{\X,p} u$, where $\lap^N_{\X,p} u=\Delta_{\X} u + (p-2)\lap_{\X,\infty} u/|\X u|^2$ is the normalized $p$-Laplacian; in the table the $p$-Laplacian is multiplied by $|\X u|^{4-p}$, which gives polynomial coefficients for every $p$. Quasi-linear equations of these types have been widely studied, see
\cite{B02, B12, C-D-G, C-D-G2, CCHP18, Cap--reg, Dom, Dom-Man--hor,Liu-Man-Zhou, Muk-Zhong, Muk0, Citti-Muk} and the references therein. On Carnot groups the viscosity solutions of the $p$-Laplace equation are its weak solutions \cite{B12}, for which the uniqueness for the Dirichlet problem is standard. Row (4) is the minimal surface equation $\dvh\big(\X u/\sqrt{1+|\X u|^2}\big)=0$, which is one of the examples of \cite{Barles-Busca}. It satisfies (E) but, for $m\geq2$, not (S): the eigenvalues of $A(\xi)$ are $1/(1+|\xi|^2)$ in the direction of $\xi$ and $1$ in the orthogonal directions, so that $A(t\xi)$ is not a multiple of $A(\xi)$ for $t>1$. Our results do not apply to it; see Remark~\ref{rem:scale}(c). Rows (5) and (7) are Aronsson equations,
\begin{equation}\label{eq:cwang}
\inp{\XXs u \,\gr f(\X u)}{\gr f(\X u)}=\sum_{i,j=1}^m \del_i f (\X u)\del_j f(\X u) X_i X_j u=0,
\end{equation}
with $f=\frac12g^2$. Equation \eqref{eq:cwang} was studied by C.~Wang \cite{CYW07} in the context of absolute minimizers of $\| f(\X u)\|_{L^\infty}$, for functions $f$ that are convex, positively homogeneous of degree $k\geq 1$ and positive off the origin; we refer to \cite{CYW07} for the precise hypotheses of the uniqueness theorem. Replacing such an $f$ by $\frac12f^{2/k}$ multiplies the operator by a positive function of $\X u$, so that there is no loss of generality in taking $k=2$. For more general supremal functionals in Carnot--Carath\'eodory spaces see \cite{Wang-yu, Pin-v-Wang, CGPV}.

Rows (6)--(10) are equations for which the approximation method is not available. Row (6) is the tilted $\infty$-Laplacian $\lap^T_{\X,\infty}u=\inp{\XXs u\;T\X u}{T\X u}$. If $T$ is symmetric this is the $\infty$-Laplacian of the frame $T^{1/2}\X$. If $T$ is not symmetric there is no function $f\in C^1(\R^m\mns\{0\})$ whose gradient is a non-vanishing multiple of $T\xi$ for all $\xi\neq0$,\footnote{Suppose that $\gr f(\xi)=\mu(\xi)T\xi$ with $\mu$ continuous and without zeros, hence of constant sign. Let $S$ be the symmetric part of $T$, $K':=\frac12S^{-1/2}(T-T^t)S^{-1/2}\neq0$, $e_1$ a unit vector with $K'e_1\neq0$ and $e_2:=K'e_1/|K'e_1|$, which is orthogonal to $e_1$. Along the closed curve $\ell(\phi)=S^{-1/2}(\cos\phi\,e_1+\sin\phi\,e_2)$ one computes $\frac{d}{d\phi}f(\ell(\phi))=\mu(\ell(\phi))\inp{(\I_m+K')c(\phi)}{c'(\phi)}=\mu(\ell(\phi))\,|K'e_1|$, where $c(\phi)=\cos\phi\,e_1+\sin\phi\,e_2$, so that $f\circ\ell$ would be strictly monotone on $[0,2\pi]$, which is impossible.} so that $\lap^T_{\X,\infty}$ is an operator of rank one that is not the Aronsson operator of any function with non-vanishing gradient. In particular it is not, even up to a positive factor, of the form \eqref{eq:cwang} for a convex function $f$ homogeneous of degree $k\ge 1$. In row (7) the function $g$ is not convex; then the functionals $\int g(\X u)^p$ are not convex and their Euler--Lagrange equations are not elliptic. Row (8) is the orthotropic $p$-Laplacian
\[
\sum_{i=1}^mX_i\big(|X_iu|^{p-2}X_iu\big)=(p-1)\sum_{i=1}^m|X_iu|^{p-2}X_iX_iu ,
\]
where we take $p\geq3$ so that its matrix is locally Lipschitz continuous; (E) and (S) hold for every $p\geq2$. It is the Euler--Lagrange equation of $\int\sum_i|X_iu|^p$, but for $m\geq2$ it is not an Aronsson operator, and its limit as $p\to\infty$ is discontinuous in $\X u$. For $p=4$ it is, up to a constant factor, the diagonal $\infty$-Laplacian
\begin{equation}\label{eq:diag}
\LL u=-\sum_{i=1}^m (X_iu)^2X_iX_iu,\qquad A(\xi)=\operatorname{diag}(\xi_1^2,\ldots,\xi_m^2).
\end{equation}
Since the orthotropic $p$-Laplacian is in divergence form and its vector field is strictly monotone, the comparison principle for its weak solutions is standard, and this is the sense of the entry in the last column of the table; for its viscosity solutions we only have Theorems~\ref{thm:comp} and \ref{thm:uncond}, as for rows (6), (7), (9) and (10). Rows (9) and (10) are the biased $\infty$-Laplacian $\lap_{\X,\infty}u+\beta|\X u|^3=0$, which in the Euclidean case is the equation of the biased tug-of-war \cite{PPS10}, and the $\infty$-Laplacian with drift $\lap_{\X,\infty}u+\inp{b}{\X u}|\X u|^2=0$, which is the $\infty$-Laplace equation of the frame $e^{\langle b,\pi_1(x)\rangle}\X$, where $\pi_1$ is the projection onto the first layer. This frame is not left-invariant, and the uniqueness of its $\infty$-harmonic functions is not covered by \cite{CYW07}. More generally, to each example one may add a Hamiltonian $H(\xi)=|\xi|\,a(\xi)H_0(\xi/|\xi|)$ of the same homogeneity along rays; see \cite{Arm-Tran} for second-order Hamilton--Jacobi equations with general Hamiltonians.

\subsection{Main results}\label{subsec:results}
Our first result is the comparison principle under the non-degeneracy hypothesis of the first step of the Barles--Busca argument, imposed on the regularizations of $u$ and $v$. For $\eps>0$ let $u^\eps$ and $v_\eps$ be the sup- and inf-convolutions of $u$ and $v$ defined in \eqref{eq:supinfconv}; they are semi-convex and semi-concave sub- and super-solutions in the reduced domains $\Om[\eps]\uparrow\Om$ of \eqref{eq:redom} (Lemma~\ref{lem:wang}). Let
\begin{equation}\label{eq:Eeps}
E_\eps:=\Big\{x\in\Om[\eps]:\ (u^\eps-v_\eps)(x)=\max_{\overline{\Om[\eps]}}(u^\eps-v_\eps)\Big\}
\end{equation}
be the set of interior maximum points of $u^\eps-v_\eps$. At the points of $E_\eps$ both $u^\eps$ and $v_\eps$ are differentiable and $\X u^\eps=\X v_\eps$ (Section~\ref{subsec:semi}).
\begin{Thm}\label{thm:comp}
Let $\Om \subset \R^n$ be a bounded open set, let $\LL$ be as in \eqref{eq:maineq} with $A_{i,j}, H\in C(\R^m\mns\{0\})$ satisfying \eqref{eq:ell} and \eqref{eq:scale}, and let $u,v\in C(\bar\Om)$ satisfy $\LL u\leq 0\leq \LL v$ in the viscosity sense in $\Om$ and $u\leq v$ on $\del \Om$. If there is a sequence $\eps_k\to0$ such that $\X u^{\eps_k}$ does not vanish at any point of $E_{\eps_k}$, then $u\leq v$ in $\Om$.
\end{Thm}
The hypothesis is the point where the Euclidean argument of \cite{Barles-Busca} and ours part company. In $\R^n$ it is not needed, because the non-degeneracy can be achieved by translating $u$ (Theorem~\ref{thm:uncond}(iii)); on a non-abelian Carnot group we do not know whether it can be removed. Section~\ref{sec:obstruction} explains why the Euclidean translation argument does not extend and formulates the open problem. Since the hypothesis bears on the regularizations of $u$ and $v$, it is not easy to check directly, except in the cases of the next theorem.
\begin{Thm}\label{thm:uncond}
Let $\Om$, $\LL$, $u$, $v$ be as in Theorem~\ref{thm:comp}, without the hypothesis on $E_{\eps_k}$. Then $u\leq v$ in $\Om$ in each of the following cases:
\begin{enumerate}[(i)]
\item $v$ is of class $C^1$ with $\X v\neq0$ in a neighbourhood of $\argmax_\Om(u-v)$, or $u$ is of class $C^1$ with $\X u\neq 0$ in a neighbourhood of $\argmax_\Om(u-v)$ (the other function being merely continuous);
\item $u$ or $v$ is constant;
\item $r=1$, i.e.\ $\G$ is the Euclidean group; more generally, $u=F\circ\pi_1$ and $v=G\circ\pi_1$ for continuous $F,G$, where $\pi_1$ is the projection onto the first layer (such functions are invariant under the derived subgroup $[\G,\G]$).
\end{enumerate}
\end{Thm}
Here $\argmax_\Om(u-v)$ is the set of points of $\Om$ at which $u-v$ attains its maximum over $\bar\Om$; when this maximum is positive it is a compact subset of $\Om$, and when it is not there is nothing to prove. Case (i) holds in particular if $v\in C^1(\Om)$ and $\X v\ne0$ in $\Om$. It is the natural replacement, in the absence of comparison with cones, of comparison with affine data: every continuous sub-solution lies below every $C^1$ super-solution without horizontal critical points, for instance below $v(x)=\inp{\lambda}{\pi_1(x)}+c$ with $\lambda\ne0$ when $H(\lambda)\geq0$. Only case (i) is obtained through the hypothesis of Theorem~\ref{thm:comp}; case (ii) follows from the strong maximum principle alone, and case (iii) from the Euclidean translation argument. Case (iii) shows that the obstruction is carried entirely by the dependence of $u$ and $v$ on the non-horizontal variables; for $r=1$ it is the comparison principle of Barles and Busca for continuous sub- and super-solutions of the equations that satisfy (E) and (S), with the Lipschitz continuity in (BB1) replaced by continuity away from the origin.

Theorems~\ref{thm:comp} and \ref{thm:uncond} apply to all the equations of Table~\ref{tab:examples} except the minimal surface equation. The solvability of the Dirichlet problem $\LL u= 0$ in $\Om$, $u=f$ on $\del\Om$, requires appropriate assumptions on the boundary, see \cite{Ar-Cr-Ju, Barles-Busca, Cr-Ish-Lions}; existence of viscosity solutions usually follows from Perron's method, and for classes of solutions to which Theorem~\ref{thm:comp} or \ref{thm:uncond} applies one obtains uniqueness.

Theorems~\ref{thm:comp} and \ref{thm:uncond}, the latter in all three cases, have been formalized in Lean~4 on top of the Mathlib library; the development is available at \cite{CompND-Lean}. The results that we quote from the literature are not re-proved there: the strong maximum principle of Bardi and Goffi in the form of Corollary~\ref{cor:stmx}, Pansu's theorem \cite{P89} and \cite[Proposition~4.8]{DMV13}, the sup- and inf-convolutions of \cite{CYW07} in the form of Lemma~\ref{lem:wang}(1)--(3), Aleksandrov's theorem and Jensen's lemma. They are stated precisely and appear as explicit hypotheses of the formalized theorems, so that each Lean theorem asserts that the cited results imply the comparison principle. Everything else is proved, and the Lean kernel confirms that no axioms beyond the standard ones of Lean are used. The repository \cite{CompND-Lean} contains the Lean sources and a blueprint with a natural-language statement and proof of every Lean declaration.

\subsection{Obstacles and key ideas}\label{subsec:ideas}
Some results about elliptic equations on the Euclidean group extend verbatim, or with routine modifications, to sub-elliptic equations on non-abelian Lie groups. This is not the case here, and the difficulty is of a specific nature.

As in \cite{Barles-Busca}, the comparison principle is approached by approximating viscosity sub/super-solutions with semi-convex/semi-concave functions, and the proof follows the two steps described in Section~\ref{subsec:BB}. The first step transfers. Under the non-degeneracy assumption that the horizontal gradient does not vanish at the maximal points of $u-v$ (Proposition~\ref{prop:comp0}), a small perturbation $u\mapsto\psi(u)$ together with the strict ellipticity \eqref{eq:ell} produces strict sub-solutions which, by Jensen's lemma and Aleksandrov's theorem, are classical sub-solutions at points of twice differentiability arbitrarily close to the maximal points. The relabelling identity behind this step (Lemma~\ref{lem:tech}) is exact under \eqref{eq:scale}, and this is why the continuity of $A$ and $H$ suffices. Theorem~\ref{thm:comp} is this step applied to the regularizations of $u$ and $v$.

The second step, the removal of the non-degeneracy assumption, is where the difficulty resides. In the Euclidean case \cite{Ar-Cr-Ju, Barles-Busca} one looks at the maximal points $x_h$ of the translates $u(\cdot+h)-v$: if the gradient vanishes at some $x_h$ for every $h$ near $0$, then $M(h)=\max(u(\cdot+h)-v)$ is constant near $0$ (Rademacher's theorem and the first-order Taylor expansion), $u$ has an interior maximum and $v$ an interior minimum, and the strong maximum principle makes both constant, contradicting the boundary inequality; otherwise there is $h$ with non-vanishing gradients at all maximal points of $u(\cdot+h)-v$, which is a sub-solution, and the first step applies to the pair $(u(\cdot+h),v)$.

On a Carnot group the two properties of the Euclidean translation used in this argument are split between two operations. Left translations $x\mapsto h*x$ preserve the equation, because the vector fields are left-invariant: if $u$, $v$ are sub/super-solutions then so are $u(h*\cdot)$, $v(h*\cdot)$ in a reduced domain. But the first-order Taylor expansion in the horizontal gradient is the one for \emph{right} translations,
$$ f(x*h)= f(x) + \inp{\X f(x)}{\pi_1(h)} + o\big(\|h\|\big),$$
where $\|x\| \sim |\pi_1(x)|+|\pi_2(x)|^{1/2}+ \ldots + |\pi_r(x)|^{1/r}$ is a homogeneous norm of $\G$. The displacement $h*x-x$ produced by a left translation is governed by the right-invariant vector fields, which the equation does not see, and it is only of order $d(h*x,x)\lesssim\|h\|^{1/r}$, so that $h\mapsto\max(u(h*\cdot)-v)$ is merely H\"older continuous when $r>1$. Right translates, on the other hand, are not sub-solutions in general: we show in Proposition~\ref{prop:central} that the horizontal vector fields are invariant under the right translation by $g$ exactly when $g$ belongs to the center of $\G$, so that $u(\cdot*g)$ is a sub-solution for every sub-solution $u$ when $g$ is central and, in general, not otherwise; and that along central directions the parameter derivative of $u(x*g)$ is a derivative of $u$ in a non-horizontal direction, so that its vanishing at maximal points gives a maximum of $u$ along a curve rather than on an open set, which the strong maximum principle does not see.

What survives of the second step is the following (Lemma~\ref{lem:case1}): for every $\eta>0$ there is a right translation $h$ with $\|h\|<\eta$ such that $\X u$ does not vanish at the points $x*h$, $x$ a maximal point of $u(\cdot*h)-v$. The proof uses the Taylor expansion above, Pansu's theorem and the strong maximum principle, as in the Euclidean case. In $\R^n$ this is all that is needed, because $u(\cdot+h)$ is a sub-solution and the first step applies to the pair $(u(\cdot+h),v)$; this is Theorem~\ref{thm:uncond}(iii). On a non-abelian group $u(\cdot*h)$ is not a sub-solution, and the first step can only be applied to $(u,v)$ itself, which requires $\X u\neq0$ on the maximal set of $u-v$. We do not know how to pass from the non-degeneracy at the maximal points of $u(\cdot*h)-v$, which is available at every scale, to the non-degeneracy at the maximal points of $u-v$: the former does not come with a lower bound as $h\to0$, and the two maximal sets are close only in a qualitative sense. Section~\ref{sec:obstruction} quantifies what is missing (Proposition~\ref{prop:transfer}), explains why the strong maximum principle alone cannot provide it, and states the open problem.

\section{Notation and preliminaries}\label{sec:prelim}
The Euclidean scalar product on $\R^n$ is $\inp{\cdot}{\cdot}$, $\gr u=(\del_{x_1}u,\ldots, \del_{x_n} u)$ is the Euclidean gradient, $DF$ the Jacobian matrix of $F:\Om \to \R^n$ and $D^2u$ the Hessian. For symmetric matrices $A, B\in \R^{k\times k}$ we write
\begin{equation}\label{eq:matord}
A\leq B \quad \text{iff} \quad \inp{A \xi}{\xi} \leq \inp{B\xi}{\xi}\quad \forall\ \xi\in \R^k;
\end{equation}
$|A|=\sqrt{\Tr(A^TA)}$ is the Frobenius norm. For non-negative definite $A, Z\geq 0$ one has $\Tr(AZ) = |\sqrt{Z}\sqrt{A}|^2\geq 0$, hence
\begin{equation}\label{eq:matmon}
\Tr(AY) \leq \Tr(AX),\qquad \forall\ A\geq 0, \ Y\leq X.
\end{equation}

\subsection{Normalization and structure conditions}\label{subsec:str}
By Remark~\ref{rem:scale}(a), under \eqref{eq:ell}--\eqref{eq:scale} we have $A(\xi)=a(\xi)A_0(\xi/|\xi|)$ and $H(\xi)=|\xi|\,a(\xi)H_0(\xi/|\xi|)$ with $a>0$ continuous, so that for $\varphi\in C^2$ and $\X\varphi(x)\ne0$
$$\LL\varphi(x)=a\big(\X\varphi(x)\big)\,\LL_0\varphi(x),\qquad \LL_0w:=-\Tr\Big(A_0\Big(\frac{\X w}{|\X w|}\Big)\XXs w\Big)+|\X w|\,H_0\Big(\frac{\X w}{|\X w|}\Big).$$
Since $a>0$ and test functions are only evaluated at points where their horizontal gradient does not vanish (Definition~\ref{def:visc} below), $\LL$ and $\LL_0$ have the same viscosity sub- and super-solutions, and $\LL_0$ satisfies \eqref{eq:ell} and \eqref{eq:scale} with $\kappa\equiv1$. Replacing $\LL$ by $\LL_0$, we assume from now on that
\begin{equation}\label{eq:normal}
A(t\xi)=A(\xi)\qquad\text{and}\qquad H(t\xi)=t\,H(\xi)\qquad\text{for all}\ t>0,\ \xi\in\R^m\mns\{0\};
\end{equation}
in particular $A$ and $H$ are bounded on bounded sets. Then $\E$ is $2$-homogeneous, $\E(t\xi)=t^2\E(\xi)$, and by \eqref{eq:ell} and the compactness of the unit sphere
\begin{equation}\label{eq:ell0}
\E(\xi)\ \geq\ e_0\,|\xi|^2\qquad\text{for all}\ \xi\in\R^m,\qquad\text{where}\quad e_0:=\min_{|\om|=1}\E(\om)>0 .
\end{equation}
Moreover, for every $t>0$, $\xi\ne0$ and every symmetric $X\in\R^{m\times m}$,
\begin{equation}\label{eq:scaleL}
-t\Tr\big(A(t\xi)X\big)+H(t\xi)= t\Big[-\Tr\big(A(\xi)X\big)+H(\xi)\Big].
\end{equation}

For smooth vector fields $X_1,\ldots,X_m$ on $\R^n$ and $u:\Om \to \R$, the horizontal gradient and the horizontal second derivative are $\X u = (X_1u, \ldots, X_m u)$ and $\XX u = (X_jX_iu)_{i,j}$, and we set
$$\XXs u =  \tfrac{1}{2}\left(\XX u +(\XX u)^T\right)= \tfrac{1}{2}\big(X_iX_ju+X_jX_iu\big)_{i,j}\,.$$
The divergence of $F= (f_1,\ldots, f_m)$ with respect to the vector fields is $\dvh (F) = \sum_{j} X_j^* f_j$, with $X_j^*$ the $L^2$-adjoint of $X_j$; for left-invariant fields in exponential coordinates $X_j^*=-X_j$. Let $\sigma:\R^n \to \R^{n\times m}$, $\sigma(x)=[\sigma^1(x),  \dots, \sigma^m(x)]$ with $\sigma^j:\R^n \to \R^n$, be such that $X_j = \sigma^j(x)\cdot \gr$. Then $\sigma$ is smooth, $\sigma(x)\neq 0$ for all $x$, and
\begin{equation}\label{eq:sigrad}
\X u= \sigma(x)^T\gr u, \qquad \XXs u = \sigma(x)^T D^2u\, \sigma(x) + \mathcal M(x,\gr u),
\end{equation}
where $\mathcal M(x,\zeta)\in \R^{m\times m}$ is the matrix with entries
\begin{equation}\label{eq:Mxz}
\mathcal M(x,\zeta)_{i,j}=\frac{1}{2}\big(\inp{D\sigma^j(x)\sigma^i(x)}{\zeta}+
\inp{D\sigma^i(x)\sigma^j(x)}{\zeta}\big),
\end{equation}
linear in $\zeta$ and smooth in $x$. On a bounded set, $|\sigma|$, $|D\sigma|$ are bounded by constants depending on $\G$ and $\diam(\Om)$; we write $c=c(\G,\Om)$ for such constants.

For $f:\Om\to \R$ we denote
\begin{equation}\label{eq:argminmax}
\argmax_\Om(f) := \big\{x\in \Om :  f(x)= \sup\nolimits_{\,\Om} f\big\}, \qquad
\argmin_\Om(f) := \big\{x\in \Om :  f(x)= \inf\nolimits_{\,\Om} f\big\};
\end{equation}
these sets may be empty. If $\Theta:\Om\to \R^n$ is a bijection onto its image, $x\in \argmax_\Om(f)$ iff $\Theta (x)\in \argmax_{\Theta(\Om)}(f\circ \inv{\Theta})$. If $f_k\to f$ uniformly, $x_k\in \argmax_\Om (f_k)$ and $x_k \to x \in \Om$, then $x\in \argmax_\Om (f)$; and $\argmax_\Om(f+c)=\argmax_\Om(f)$ for constants $c$.

\subsection{Viscosity solutions of sub-elliptic equations}\label{subsec:nonsubeq}
For the classical theory of viscosity solutions we refer to \cite{Ar-Cr-Ju, Cr-Ish-Lions, J88}. A sub-elliptic equation
\begin{equation}\label{eq:G}
G(x,u,\X u, \XXs u)=0,\qquad G: \bar\Om \times \R \times (\R^m\mns\{0\})\times \R^{m\times m}\to \R,
\end{equation}
can be written as $F(x,u, \gr u, D^2 u)=0$ using \eqref{eq:sigrad}, with
\begin{equation}\label{eq:FtoG}
F(x,r, \zeta, X):= G\big(x,r,\, \sigma(x)^T\zeta,\, \sigma(x)^T X \sigma(x) + \mathcal M(x,\zeta)\big),
\end{equation}
so that the notion of viscosity solution can be defined as in the Euclidean case; since $G$ may be singular or undefined where the horizontal gradient vanishes, as is the case for the normalized $p$-Laplacian, test functions are restricted accordingly. Let $\USC$ and $\LSC$ denote the classes of upper and lower semi-continuous functions. For $w:\Om\to\R$ and $x\in \Om$, let
\begin{equation}\label{eq:defA}
\Av_x^+ (w, \Om) =\{\varphi\in C^2(\Om) : x\in \argmax_\Om(w-\varphi),\ \X \varphi(x)\neq 0\},
\end{equation}
and let $\Av_x^-(w,\Om)$ be defined similarly with $\argmin$ in place of $\argmax$. Note that $\X\varphi(x)=\sigma(x)^T\gr\varphi(x)\neq0$ implies $\gr\varphi(x)\neq0$; the converse fails when $m<n$.
\begin{Def}\label{def:visc}
For the equation \eqref{eq:G}, $u\in \USC(\Om)$ (resp. $u\in \LSC(\Om)$) is a viscosity sub-solution (resp. super-solution) in $\Om$ if for every $x\in\Om$ and every $\varphi\in \Av_x^+(u,\Om)$ (resp. $\Av_x^-(u,\Om)$),
$$ G(x,\varphi(x),\X \varphi(x), \XXs \varphi(x)) \leq 0 \ (\text{resp.} \geq 0);$$
we write $G(x,u,\X u, \XXs u)\leq 0$ (resp. $\geq 0$) in the viscosity sense. A viscosity solution is both a sub- and a super-solution.
\end{Def}
Since $\varphi$ may be modified away from $x$, it is equivalent to require in \eqref{eq:defA} that $x$ be a local maximum of $w-\varphi$. For the operator \eqref{eq:maineq}, $u$ is a sub-solution if $\LL\varphi(x)\leq0$ for all $\varphi\in\Av^+_x(u,\Om)$.
\begin{Rem}\label{rem:classical}
If $u$ is a viscosity sub-solution (resp. super-solution) of \eqref{eq:maineq}, twice differentiable at $x$ with $\X u(x)\ne0$, then $\LL u(x)\le0$ (resp. $\geq0$) classically. Indeed, for $\eta>0$ the second-order Taylor polynomial of $u$ at $x$ plus $\eta|y-x|^2$ (resp. minus) belongs to $\Av^+_x(u,\Om)$ (resp. $\Av^-_x$) after a modification away from $x$, and $\LL$ is continuous in $(\X\varphi,\XXs\varphi)$ at points with $\X\varphi\ne0$; let $\eta\to0^+$.
\end{Rem}

The main external result we rely on is the strong maximum principle of Bardi and Goffi \cite{Bardi-Goffi} for fully non-linear sub-elliptic equations. We state it in the form we need, for sub- and super-solutions.
\begin{Thm}[Strong maximum and minimum principles, \cite{Bardi-Goffi}]\label{thm:strongmax}
Let $X_1,\ldots,X_m$ be smooth vector fields satisfying H\"ormander's condition \eqref{eq:horcond}, and let $G: \bar\Om \times \R \times (\R^m\mns\{0\})\times \R^{m\times m}\to \R$ satisfy:
\begin{enumerate}
\item $G$ is lower semicontinuous and proper: $G(x,r,\xi, X)\leq G(x,s,\xi, Y)$ for all $r\leq s$ and symmetric $Y\leq X$;
\item (scaling) for all $x\in\Om$, $\xi\ne0$, symmetric $X$, $r\in[-1,0]$ and $\lambda\in(0,1]$, the inequality $G(x,r,\xi,X)>0$ implies $G(x,\lambda r,\lambda \xi, \lambda X)>0$;
\item (ellipticity in the gradient direction) $\sup_{\gamma>0} G\,\big(x,0, \xi, X-\gamma \xi\otimes\xi\big) >0$ for all $x\in \Om$, $\xi\ne0$ and symmetric $X$.
\end{enumerate}
Then any viscosity sub-solution of $G(x,u,\X u, \XXs u)=0$ in $\Om$ that attains a non-negative maximum in $\Om$ is constant. If $G$ is upper semicontinuous and proper, satisfies (3), and satisfies (2) with ``$<0$'' in place of ``$>0$'' and $r\in[0,1]$, then any viscosity super-solution that attains a non-positive minimum in $\Om$ is constant.
\end{Thm}
In \cite{Bardi-Goffi} the scaling hypothesis is written as $G(x,\lambda r,\lambda \xi, \lambda X) \geq \phi(\lambda)G(x,r,\xi, X)$ for some $\phi:(0,1]\to(0,\infty)$; by \cite[Remark 2.8]{Bardi-Goffi} the function $\phi$ may depend on $(x,r,\xi,X)$, and only the sign preservation in (2) is used; the super-solution version is \cite[\S2.3]{Bardi-Goffi}. Viscosity solutions are defined in \cite{Bardi-Goffi} with test functions $\varphi$ satisfying $\gr\varphi(x)\ne0$, since their $G$ is defined for all non-zero Euclidean gradients, whereas our $F$ in \eqref{eq:FtoG} is defined only for $\sigma(x)^T\zeta\ne0$ and our class of test functions \eqref{eq:defA} is smaller, so that our sub-solutions form a larger class. We therefore check that the proof of \cite[Proposition 2.4]{Bardi-Goffi}, from which the strong maximum principle follows by Chow's theorem, uses only test functions of our class. In that proof $K$ is the set where the sub-solution $u$ attains its maximum, $z\in\del K\cap\Om$ is a point at which $K$ has an exterior Euclidean ball $B(y,R)$, $\nu:=(y-z)/R$, and the sub-solution property is used only for the test functions $\eps v$, $v(x)=e^{-\gamma R^2}-e^{-\gamma|x-y|^2}$, at points of a small Euclidean ball $B(z,\rho)$: at a maximum point of $u-\eps v$ in $B(z,\rho)\cap B(y,R)$ and at $z$ itself. Since $\gr v(x)=2\gamma e^{-\gamma|x-y|^2}(x-y)$, the horizontal gradient of $v$ at $z$ is $-2\gamma Re^{-\gamma R^2}\sigma(z)^T\nu$, and the hypothesis to be contradicted in that proof is precisely $\sigma^i(z)\cdot\nu\ne0$ for some $i$, that is, $\sigma(z)^T\nu\ne0$ (the columns $\sigma^i$ are the subunit vector fields of the equation, \cite[Lemma 3.1]{Bardi-Goffi}). Hence $\X v\ne0$ at $z$ and, by continuity, on $B(z,\rho)$ for $\rho$ small, so that $F$ is defined and continuous there, all the test functions used belong to \eqref{eq:defA}, and the proof applies verbatim.
\begin{Cor}\label{cor:stmx}
Let $\LL$ be as in \eqref{eq:maineq} with $A$, $H$ satisfying \eqref{eq:ell}--\eqref{eq:scale}, and let $\Om$ be connected. Any viscosity sub-solution of $\LL u=0$ in $\Om$ that attains a non-negative maximum in $\Om$ is constant, and any viscosity super-solution that attains a non-positive minimum in $\Om$ is constant.
\end{Cor}
\begin{proof}
$G(x,r,\xi, X) = -\Tr \big(A(\xi) X\big)+H(\xi)$ is continuous on its domain and proper by \eqref{eq:matmon}, since $A(\xi)\geq0$. By the normalization \eqref{eq:normal},
$$G(x,\lambda r,\lambda\xi,\lambda X)=-\lambda\Tr\big(A(\lambda \xi) X\big) + H(\lambda \xi)
 = \lambda\,G(x,r,\xi,X),$$
so that both sign conditions in (2) hold. (Without the normalization, \eqref{eq:scale} gives the positive factor $1/(t\kappa(t,\lambda\xi))$ with $t=1/\lambda$ instead of $\lambda$.) Finally
$$  G\big(x,0, \xi, X-\gamma \xi\otimes\xi\big)
\,=\, \gamma \E(\xi) -\Tr \big(A(\xi) X\big) + H(\xi)>0  $$
for $\gamma$ large, by \eqref{eq:ell}. Theorem~\ref{thm:strongmax} applies.
\end{proof}

\subsection{Semi-concave and semi-convex functions}\label{subsec:semi}
In this subsection balls and distances are Euclidean.
\begin{Def}\label{def:semicon}
A function $w\in C(\Om)$ is semi-convex in $\Om$ if there exists $\Lambda\geq0$ such that $x\mapsto w(x) + \frac{1}{2}\Lambda |x|^2$ is convex on every Euclidean ball $B\subset\Om$; $w$ is semi-concave if $-w$ is semi-convex. The smallest such $\Lambda$ is the semi-convexity constant of $w$.
\end{Def}
We shall use the following properties, see \cite{Barles-Busca, Cr-Ish-Lions, Evans-Gar}.
\begin{enumerate}
\item (Differentiability at maximal points) If $u$ is semi-convex and $v$ semi-concave in $\Om$, then $u$ and $v$ are differentiable at every point of $\argmax_\Om(u-v)$, with $\gr u=\gr v$ there, and at every local maximum point of $u-v$.
\item (Partial continuity of the gradient) If $w$ is semi-convex or semi-concave, differentiable at $x$ and at the points of a sequence $x_k\to x$, then $\gr w(x_k) \to \gr w(x)$.
\item (Local Lipschitz bound) If $w$ is semi-convex in $\Om$ with constant $\Lambda$ and $K\subset\Om$ is compact, then $w$ is Lipschitz on $K$, with constant bounded by $c(n)\,\big(\osc_\Om w+\Lambda\sup_\Om|x|\diam(\Om)\big)\,\dist(K,\del\Om)^{-1}$; in particular $\gr w$ is bounded on $K$ where it exists.
\item (Composition) If $w$ is semi-convex, $\psi\in C^2(\R)$ with $\psi'\ge1$, $\psi''\geq0$, and $\Theta$ is a $C^2$ diffeomorphism with bounded first and second derivatives, then $\psi\circ w$ and $w\circ\Theta$ are semi-convex on compact subsets, with constants depending on $\Lambda$, on the local Lipschitz constant of $w$, and on $\psi$, $\Theta$.
\end{enumerate}
Property (3) is the standard Lipschitz bound for the convex function $w+\frac12\Lambda|x|^2$; (4) follows from $D^2(\psi\circ w)=\psi''\gr w\otimes\gr w+\psi'D^2w$ and $D^2(w\circ\Theta)=D\Theta^TD^2w\,D\Theta+\gr w\cdot D^2\Theta$ at points of twice differentiability, and by approximation elsewhere.
\begin{Thm}[Aleksandrov]\label{thm:aleksandrov}
A semi-convex function is twice differentiable almost everywhere.
\end{Thm}
The set $\argmax(w)$ may have measure zero, so twice differentiability at maximal points is not guaranteed. The following lemma (\cite[Lemma A.3]{Cr-Ish-Lions}) provides maximal points of linear perturbations, arbitrarily close to a strict maximum point, at which the second derivatives exist.
\begin{Lem}[Jensen's Lemma]\label{lem:jensen}
Let $w$ be semi-convex in a neighbourhood of $\hat x$, and let $\hat x$ be a strict local maximum point of $w$. For $p\in \R^n$ let $w_p(x)= w(x) + \inp{p}{x}$. Then, for all $\rho,\delta>0$ small enough, the set
\begin{equation}\label{eq:krdelta}
J_{\rho,\delta}(\hat x)= \bigcup_{|p|< \delta} B_\rho(\hat x)\cap\argmax_{B_\rho(\hat x)}(w_p)
\end{equation}
has positive Lebesgue measure.
\end{Lem}
Since $D^2 w_p = D^2 w$ wherever the latter exists, Aleksandrov's theorem gives, for any $\rho,\delta>0$ small, a point $z$ with $|z-\hat x|<\rho$ and $p$ with $|p|<\delta$ such that $w$ is twice differentiable at $z$ and $z$ is a maximum point of $w_p$ over $B_\rho(\hat x)$. The strictness of the maximum cannot be dispensed with (if $w$ is constant near $\hat x$, $w_p$ has no interior maximum in $B_\rho(\hat x)$); in Proposition~\ref{prop:comp0} it is obtained by a quartic penalization.

\subsection{Carnot groups}\label{subsec:carnot}
We refer to \cite{Bonfig-Lanco-Ugu} for the structure and properties of homogeneous and stratified groups.
\begin{Def}\label{def:carnot}
A connected, simply connected Lie group $\G=(\R^n, *)$ is a Carnot group (stratified group) of step $r$ with $m$ generators if its Lie algebra $\g$ admits a decomposition $\g= \g_1\oplus \ldots \oplus \g_r$ with $\dim\g_1=m$ and
\begin{equation}\label{eq:strat}
[\,\g_1, \g_{j}\,] = \g_{j+1}\quad\text{for}\ j\in\{1,\ldots,r-1\},\qquad [\,\g_1,\g_r\,] = \{0\}.
\end{equation}
\end{Def}
The stratification determines the dilations $\triangle_\lambda \big(Y_1+\ldots+Y_r\big)= \lambda Y_1+\ldots+\lambda^rY_r$ ($Y_j\in\g_j$), a one-parameter family of automorphisms of $\g$, and $\delta_{\lambda} = \exp\,\circ\,\triangle_{\lambda}\circ\exp^{-1}$ on $\G$, where $\exp: \g \to \G$ is a global diffeomorphism. Conversely, if $\G=(\R^n,*)$ is a homogeneous Lie group and $X_1,\ldots,X_m$ are linearly independent left-invariant vector fields, $1$-homogeneous with respect to the dilations, satisfying \eqref{eq:horcond}, then with $\g_1=\spn\{X_1,\ldots,X_m\}$ and $\g_j$ the span of the $j$-fold brackets of elements of $\g_1$, the sum $\g_1+\ldots+\g_r$ is direct (the $\g_j$ lie in distinct eigenspaces of the dilations), it equals $\g$ by \eqref{eq:horcond}, and \eqref{eq:strat} holds; so $\G$ is a Carnot group and $X_1,\ldots,X_m$ a basis of its first layer. (Homogeneity is essential here: there are nilpotent Lie algebras generated by $m$ elements that admit no stratification, e.g.\ the $5$-dimensional algebra with basis $x_1,\ldots,x_5$ and non-zero brackets $[x_1,x_2]=x_3$, $[x_1,x_3]=x_4$, $[x_1,x_4]=x_5$, $[x_2,x_3]=x_5$: its lower central series has quotients of dimensions $2,1,1,1$, while $[\g_1,[\g_1,\g_1]]$ is $2$-dimensional for every complement $\g_1$ of the derived algebra.)

Letting $m_1=m$ and $m_j=\dim(\g_j)$, we extend $X_1,\ldots,X_m$ to a basis $\{X_{i,j}\}$ of $\g$ with $X_i=X_{i,1}$ and $\g_j=\spn\{ X_{i,j} : 1\leq i\leq m_j\}$, orthonormal for a chosen Euclidean norm on $\g$ for which the $\g_j$ are orthogonal. In exponential coordinates $x= \exp\big( \sum_{j=1}^r \sum_{i=1}^{m_j}  x_{i,j} X_{i,j}\big)$, the group law $(x,y)\mapsto x*y$ is polynomial, $\inv{x}=-x$, $X_i=\del_{x_{i,1}}+\sum_{j\geq2}\sum_k b^i_{k,j}(x)\del_{x_{k,j}}$ with polynomial $b^i_{k,j}$ depending only on coordinates of layers $<j$, and
\begin{equation}\label{eq:dilG}
\delta_\lambda (x)= (\lambda x_{1}, \lambda^2  x_{2}, \ldots, \lambda^r x_{r}) \quad\text{where}\quad  x_j= ( x_{1,j},\ldots,  x_{m_j,j})\in \R^{m_j}.
\end{equation}
We denote by $\pi_j:\R^n\to \R^{m_j}$ the projection $x\mapsto x_j$; $\pi_1(x*y)=\pi_1(x)+\pi_1(y)$, so $\pi_1$ is a homomorphism onto the abelianization $\G/[\G,\G]\cong\R^m$, and $X_iu(x) =\lim_{t\to 0} \frac{1}{t}\big(u(x*\exp(tX_i))-u(x)\big)$ for $u\in C^1$. We fix the homogeneous norm
\begin{equation}\label{eq:homnorm}
\|x\| = \bigg(\sum_{j=1}^r \Big(\sum_{i=1}^{m_j}| x_{i,j}|^2\Big)^\frac{r!}{j}
\bigg)^\frac{1}{2r!},
\end{equation}
which is continuous, smooth off the origin, even, and satisfies $\|\delta_\lambda x\|=\lambda\|x\|$. It induces the left-invariant metric $d(x,y)=\|\inv{y}* x\|$, which satisfies $d (z*x,z*y) = d (x,y) $, $d (\delta_\lambda x,\delta_\lambda y) = \lambda d (x,y)$, and
\begin{equation}\label{eq:righttrans}
d(x*h,x)=\|h\|\qquad\text{for all}\ x,h\in\G .
\end{equation}
We write $B_\rho(x)=\{y: d(y,x)<\rho\}$ and $\dist(x,E)=\inf\{d(x,y): y\in E\}$; the balls $B_\rho(0)$ are connected (they are star-shaped with respect to the dilations), hence so are all balls. On bounded sets $d$ and the Euclidean distance are H\"older equivalent: $c^{-1}|x-y|\le d(x,y)\le c|x-y|^{1/r}$ for $x,y$ in a bounded set, and the topology is the Euclidean one. If $\|x\|,\|y\|<\nu$, there exists $c_\G=c(\G,\nu)>0$ such that
\begin{equation}\label{eq:trinangle}
\|x*y\| \leq c_\G\, (\|x\|+ \|y\|),\qquad\text{and}\qquad
\|\inv{x}*y*x\| \leq c_\G\, \|y\|^\frac{1}{r},
\end{equation}
see \cite{Bonfig-Lanco-Ugu, Mag-strat}; the second inequality in \eqref{eq:trinangle} is the estimate of the inner automorphism $\cc_x(y)=\inv{x}*y*x$, and the exponent $1/r$ is sharp. For $u$ differentiable at $x_0$ we have the first-order Taylor formula
\begin{equation}\label{eq:taylor2}
u(x_0*h)= u(x_0) + \inp{\X u(x_0)}{\pi_1(h)} + o\big(\|h\|\big)\qquad(h\to0),
\end{equation}
which follows from the Euclidean expansion of $u$ at $x_0$ and from the identity $x_0*h-x_0=\sum_{i=1}^m\pi_1(h)_i\,X_i(x_0)+O(\|h\|^2)$: to first order, the displacement produced by the right translation $R_h$ is the horizontal vector $\sum_i\pi_1(h)_iX_i$ evaluated at $x_0$. (For left translations the displacement $h*x_0-x_0$ involves the right-invariant vector fields instead, and this is why $u(h*x_0)$ has no expansion in terms of $\X u(x_0)$; see Section~\ref{subsec:ideas}.) By Pansu's theorem \cite{P89}, a function $f$ that is Lipschitz with respect to $d$ on an open set is Pansu differentiable at almost every point $h$: there is a vector $\X f(h)\in\R^m$, the horizontal differential of $f$ at $h$, such that $f(h*g)=f(h)+\inp{\X f(h)}{\pi_1(g)}+o(\|g\|)$ as $g\to0$ (for $f\in C^1$ this is \eqref{eq:taylor2}); and a $d$-Lipschitz function on a connected open set whose horizontal differential vanishes almost everywhere is constant (\cite[Proposition 4.8]{DMV13}).

The Haar measure of $\G$ is the Lebesgue measure of $\R^n$, and $|\delta_\lambda (E)|=\lambda^Q|E|$ where $Q=\sum_{j=1}^r jm_j$ is the homogeneous dimension.

We record how the two translations interact with the equation. For $h\in\G$ let $L_h(x)=h*x$ and $R_h(x)=x*h$, and for $u:\Om\to\R$ let
$$u^h:=u\circ L_h\quad\text{on}\ \inv{L_h}(\Om)=\inv{h}*\Om,\qquad u_h:=u\circ R_h\quad\text{on}\ \inv{R_h}(\Om)=\Om*\inv{h}.$$
\begin{Lem}\label{lem:lefttrans}
Let $u$ be a viscosity sub-solution (resp. super-solution) of \eqref{eq:maineq} in $\Om$. Then $u^h$ is a viscosity sub-solution (resp. super-solution) in $\inv{h}*\Om$, for every $h\in\G$.
\end{Lem}
\begin{proof}
Since the $X_i$ are left-invariant, $X_i(\varphi\circ L_h)=(X_i\varphi)\circ L_h$ for $\varphi\in C^2$, hence $\LL(\varphi\circ L_h)=(\LL\varphi)\circ L_h$. If $\bar\varphi\in\Av^+_{\bar y}(u^h,\inv{h}*\Om)$, then $\varphi:=\bar\varphi\circ\inv{L_h}\in\Av^+_{h*\bar y}(u,\Om)$ ($\X\varphi(h*\bar y)=\X\bar\varphi(\bar y)\ne0$ and $u-\varphi=(u^h-\bar\varphi)\circ\inv{L_h}$ attains its maximum at $h*\bar y$), so $\LL\bar\varphi(\bar y)=\LL\varphi(h*\bar y)\le0$.
\end{proof}
Right translates $u_h$ are not sub-solutions in general; see Proposition~\ref{prop:central}.

\section{Comparison principles}\label{sec:comp}
In this section we prove comparison principles for viscosity sub/super-solutions of
\begin{equation}\label{eq:defL}
\LL u =-\Tr\big(A(\X u) \XXs u\big)+ H(\X u)=0\qquad\text{in}\ \Om,
\end{equation}
with $A$ and $H$ as in \eqref{eq:ell}--\eqref{eq:scale}, normalized as in \eqref{eq:normal}, leading to the proofs of Theorems \ref{thm:comp} and \ref{thm:uncond}. Throughout this section, $\Om\subset\R^n$ is a bounded open set, not necessarily connected.

Since $A, H\in C(\R^m\mns\{0\})$, for every compact $\mathcal Q\subset \R^m\mns\{0\}$ the local moduli of continuity
\begin{equation}\label{eq:modHA}
\begin{aligned}
\om_A(t)&:=\sup\{ |A(\xi)-A(\zeta)|\, :\, \xi,\zeta\in \mathcal Q,\ |\xi-\zeta|\leq t \},\\
\om_H(t)&:=\sup\{ |H(\xi)-H(\zeta)|\, :\, \xi,\zeta\in \mathcal Q,\ |\xi-\zeta|\leq t \},
\end{aligned}
\end{equation}
are non-decreasing and tend to $0$ as $t\to0^+$. In the arguments below the gradient variable stays in a compact set of the form $\mathcal Q=\{\theta'\le|\xi|\le C\}$ with $\theta'>0$, and we use \eqref{eq:modHA} for that set.

\subsection{Comparison for semi-convex/semi-concave functions}\label{subsec:semicomp}
In this subsection $u$ is semi-convex and $v$ semi-concave in $\Om$, both in $C(\bar\Om)$, so that they are differentiable at the maximal points of $u-v$, with $\gr u=\gr v$ and hence $\X u=\X v$ there, and their gradients are partially continuous (Section~\ref{subsec:semi}). The elementary observation on which everything rests is the following. If $u-v$ has a local maximum at a point $x_0\in\Om$ at which $u$ and $v$ are twice differentiable, then $\gr u(x_0)=\gr v(x_0)$ and $D^2u(x_0)\leq D^2v(x_0)$, hence by \eqref{eq:sigrad} $\X u(x_0)=\X v(x_0)=:\xi_0$ and $\XXs u(x_0)\leq\XXs v(x_0)$; if $\xi_0\ne0$, then \eqref{eq:matmon} gives
$$0\, \leq \Tr\Big(A(\xi_0) \big(\XXs v(x_0)-\XXs u(x_0)\big)\Big)= \LL u(x_0)-\LL v (x_0),$$
since the terms $H(\xi_0)$ cancel, so that $\LL u(x_0)<0\leq\LL v(x_0)$ is impossible. The strategy is therefore to construct perturbations $\psi(u)$ of a sub-solution $u$ which are strict sub-solutions at points of twice differentiability where the horizontal gradient does not vanish, and to produce such points near the maximal points of $\psi(u)-v$ by Jensen's lemma.

\subsubsection{The non-degenerate case}
\begin{Lem}\label{lem:tech}
Let $\psi\in C^2(\R)$ with $\psi'>0$, and let $w:\Om\to \R$ be twice differentiable at $x_0\in \Om$ with $\X w(x_0)\neq 0$. Then, with $t=\psi'(w(x_0))$ and $\xi_0=\X w(x_0)$,
\begin{equation}\label{eq:estLh}
 \LL (\psi\circ w)(x_0) = t\Big[ \LL w(x_0) -\frac{\psi''(w(x_0))}{\psi'(w(x_0))} \E\big(\xi_0\big)\Big].
\end{equation}
\end{Lem}
\begin{proof}
At $x_0$ we have $\X(\psi\circ w)=\psi'(w)\X w=t\xi_0$ and $\XXs (\psi\circ w)= \psi'(w)\XXs w + \psi''(w) \X w \otimes \X w$. Hence, by \eqref{eq:normal},
\begin{align*}
\LL(\psi\circ w)(x_0)&= -t\Tr\big(A(t\xi_0) \XXs w\big) - \psi''(w)\inp{A(t\xi_0)\xi_0}{\xi_0} +H( t\xi_0)\\
&= t\Big[ -\Tr\big(A(\xi_0) \XXs w\big) + H(\xi_0)\Big]-\psi''(w)\E(\xi_0),
\end{align*}
which is \eqref{eq:estLh}.
\end{proof}
\begin{Cor}\label{cor:gamma}
Let $\psi\in C^2(\R)$ with $\psi'> 0$ and $\psi''\geq \lambda> 0$, let $\theta>0$, and let $w:\Om\to \R$ be twice differentiable at $x_0\in \Om$ with $|\X w(x_0)|\geq \theta$ and $\LL w(x_0)\leq 0$. Then
$$ \LL (\psi\circ w)(x_0)\leq -\gamma,\qquad \gamma:=\lambda\, e_0\,\theta^2>0,$$
with $e_0$ as in \eqref{eq:ell0}.
\end{Cor}
\begin{proof}
By \eqref{eq:estLh} and $\LL w(x_0)\le0$, $\LL(\psi\circ w)(x_0)\le -\psi''(w(x_0))\,\E(\xi_0)\le-\lambda e_0|\xi_0|^2\le-\lambda e_0\theta^2$.
\end{proof}
Lemma~\ref{lem:tech} and Corollary~\ref{cor:gamma} have obvious super-solution versions, with $\psi''\le-\lambda<0$, which we shall not need: the hypothesis of Proposition~\ref{prop:comp0} below is symmetric in $u$ and $v$, because $\X u=\X v$ at the maximal points of $u-v$.

\begin{Prop}[non-degenerate case]\label{prop:comp0}
Let $u, v\in C(\bar \Om)$ be respectively semi-convex and semi-concave in $\Om$, with $u\leq v$ on $\del\Om$ and $\LL u\leq 0\leq \LL v$ in $\Om$ in the viscosity sense. If $\X u$ does not vanish at any point of $\argmax_\Om(u-v)$, then $u\leq v$ in $\Om$.
\end{Prop}
\begin{proof}
Suppose, by contradiction, that $M_0:=\max_{\bar\Om}(u-v)>0$; then $E:=\argmax_\Om(u-v)$ is a non-empty compact subset of $\Om$, since $u-v\le0$ on $\del\Om$. Replacing $v$ by $v+\tau$ with $0<\tau<M_0$ we may assume $u-v\le-\tau$ on $\del\Om$ and $M_0-\tau>0$; this does not change $E$, $\X u$, $\X v$.

\emph{Step 1: relabelling.} Since $E$ is compact, $u$ is differentiable on $E$ and $\X u$ is continuous on $E$ by partial continuity, $\theta:=\min_E|\X u|>0$. For $\lambda>0$ let
$$\psi_\lambda (s)= s+\lambda (s-u_0)^2,\qquad u_0=\inf_\Om u,\qquad u_\lambda:=\psi_\lambda\circ u ,$$
so that $\psi_\lambda'\ge1$ on $[u_0,\infty)$, $\psi_\lambda''=2\lambda$, $\|u_\lambda-u\|_{L^\infty}\le4\lambda\|u\|^2_{L^\infty}$, and $u_\lambda$ is semi-convex (Section~\ref{subsec:semi}(4)). For $\lambda$ small, $u_\lambda-v\le-\tau/2$ on $\del\Om$ and $\max_\Om(u_\lambda-v)\ge M_0-\tau-4\lambda\|u\|^2_\infty>0$, so $E_\lambda:=\argmax_\Om(u_\lambda-v)$ is a non-empty compact subset of $\Om$; moreover $E_\lambda\to E$ in the sense that $\sup_{E_\lambda}\dist(\cdot,E)\to0$ as $\lambda\to0^+$, because $u_\lambda-v\to u-v$ uniformly. Every $x\in E_\lambda$ is a differentiability point of $u_\lambda$, hence of $u$. We claim that for $\lambda$ small,
\begin{equation}\label{eq:thetahalf}
|\X u(x)|\ge\theta/2\qquad\text{for all}\ x\in E_\lambda .
\end{equation}
Otherwise there are $\lambda_k\to0$, $x_k\in E_{\lambda_k}$ with $|\X u(x_k)|<\theta/2$; a subsequence converges to some $\hat x\in E$, and partial continuity of the gradient at $\hat x$ gives $\X u(x_k)\to\X u(\hat x)$, $|\X u(\hat x)|\ge\theta$, a contradiction. Fix such a $\lambda$, and write $\tilde u=u_\lambda$, $\tilde M:=\max_\Om(\tilde u-v)>0$.

\emph{Step 2: penalization.} Fix $x_1\in E_\lambda$, let $P(x)=\mu|x-x_1|^4$ with $\mu>0$, and $\tilde v:=v+P$, which is semi-concave. Then $x_1$ is the unique maximum point of $\tilde u-\tilde v$ on $\bar\Om$, since $(\tilde u-\tilde v)(x)\le\tilde M-\mu|x-x_1|^4$.

\emph{Step 3: Jensen's lemma.} By Lemma~\ref{lem:jensen} and Theorem~\ref{thm:aleksandrov} applied to $w=\tilde u-\tilde v$ at the strict maximum $x_1$, for every $k\in\N$ large there are $p_k\in\R^n$ with $|p_k|\le1/k$ and $z_k\in\Om$ with $|z_k-x_1|<1/k$ such that $\tilde u$ and $v$ are twice differentiable at $z_k$ and $z_k$ is a local maximum point of $\tilde u-\tilde v+\inp{p_k}{\cdot}$. Hence
\begin{equation}\label{eq:zkcond}
\gr \tilde u(z_k)=\gr v(z_k)+\gr P(z_k)-p_k,\qquad D^2\tilde u(z_k)\le D^2v(z_k)+D^2P(z_k),
\end{equation}
with $|\gr P(z_k)|\le4\mu k^{-3}$, $|D^2P(z_k)|\le12\mu k^{-2}$. By \eqref{eq:sigrad}, setting $\xi_k:=\X v(z_k)$ and $\eta_k:=\sigma(z_k)^T(\gr P(z_k)-p_k)$, $|\eta_k|\le c/k$,
\begin{equation}\label{eq:klam}
\begin{aligned}
&\X \tilde u(z_{k})=\xi_k+\eta_k,\qquad \XXs \tilde u(z_{k})\leq \XXs v(z_{k})+N_k,\\
&N_k:=\sigma(z_k)^TD^2P(z_k)\,\sigma(z_k)+\mathcal M\big(z_k,\gr P(z_k)-p_k\big),\qquad |N_k|\le c/k .
\end{aligned}
\end{equation}
Since $z_k\to x_1$, $x_1$ is a differentiability point of $u$ and the $z_k$ are differentiability points of $u$, partial continuity and \eqref{eq:thetahalf} give $|\X u(z_k)|\ge\theta/4$ for $k$ large, hence $|\X\tilde u(z_k)|=\psi_\lambda'(u(z_k))|\X u(z_k)|\ge\theta/4$ and $|\xi_k|\ge\theta/4-c/k\ge\theta/8$. In particular $\X u(z_k)\ne0$ and $\X v(z_k)\ne0$, so by Remark~\ref{rem:classical} the viscosity inequalities hold classically at $z_k$: $\LL u(z_k)\le0\le\LL v(z_k)$, and by Corollary~\ref{cor:gamma} (with $w=u$, $\psi=\psi_\lambda$, $\theta/4$ in place of $\theta$)
\begin{equation}\label{eq:gamsoln}
\LL \tilde u (z_k)\leq -\gamma<0,\qquad \gamma=2\lambda\,e_0\,(\theta/4)^2 .
\end{equation}

\emph{Step 4: conclusion.} From semi-convexity of $\tilde u$ and semi-concavity of $v$ and \eqref{eq:zkcond}, $-\Lambda'\I\le D^2\tilde u(z_k)\le D^2v(z_k)+D^2P(z_k)\le(\Lambda+1)\I$ for $k$ large, so $|\XXs\tilde u(z_k)|\le C$ with $C$ independent of $k$. Using \eqref{eq:matmon} with the non-negative matrix $A(\xi_k)$ and \eqref{eq:klam},
\begin{align*}
0\, &\leq \Tr\Big(A(\xi_{k}) \big(\XXs v(z_{k})+N_k-\XXs \tilde u(z_{k})\big)\Big)\\
&= \big[-\LL v(z_k)+H(\xi_k)\big]-\Tr\big(A(\xi_k+\eta_k)\XXs\tilde u(z_k)\big)\\
&\qquad+\Tr\big([A(\xi_k+\eta_k)-A(\xi_k)]\XXs\tilde u(z_k)\big)+\Tr\big(A(\xi_k)N_k\big)\\
&= -\LL v(z_k)+\LL\tilde u(z_k)+\big[H(\xi_k)-H(\xi_k+\eta_k)\big]\\
&\qquad+\Tr\big([A(\xi_k+\eta_k)-A(\xi_k)]\XXs\tilde u(z_k)\big)+\Tr\big(A(\xi_k)N_k\big)\\
&\le -\gamma+\om_H(c/k)+C\,\om_A(c/k)+c\,\|A\|_{L^\infty(\mathcal Q)}/k ,
\end{align*}
where $\mathcal Q=\{\theta/8\le|\xi|\le C'\}$ contains all $\xi_k,\xi_k+\eta_k$ (the gradients are bounded by Section~\ref{subsec:semi}(3)). For $k$ large the right-hand side is negative, a contradiction.
\end{proof}

\subsubsection{The degenerate case: right translates}\label{subsubsec:degenerate}
To remove the non-degeneracy assumption of Proposition~\ref{prop:comp0} we study how the maximal points of $u-v$ move under right translations of $u$. For $\delta\geq0$ let
\begin{equation}\label{eq:omdelta}
\Om_\delta:=\{x\in \Om: \dist(x,\del\Om)>\delta\},\qquad K_\rho:=\{x\in\Om:\dist(x,\del\Om)\geq\rho\}\ \ (\rho>0),
\end{equation}
so that $K_\rho$ is compact and $\Om_\delta\supset K_\rho$ for $\rho>\delta$. By \eqref{eq:righttrans}, $x*h\in\Om$ whenever $x\in\Om_\delta$ and $\|h\|<\delta$. For such $h$ we consider the right translate $u_h=u(\cdot*h)$ on $\Om_\delta$ and
\begin{equation}\label{eq:MAl}
M_\delta(h):= \sup_{x\in \Om_\delta} \, \{u_h(x)-v(x)\},\qquad
\A_\delta(h):= \{x\in \Om_\delta: u_h(x)-v(x)= M_\delta(h)\}.
\end{equation}
The function $u_h$ is semi-convex on compact subsets of $\Om_\delta$ (Section~\ref{subsec:semi}(4)), so at every point $x$ of $\A_\delta(h)$ the functions $u_h$ and $v$ are differentiable, and $u$ is differentiable at $x*h$.

\emph{Standing assumptions and notation for the rest of Section~\ref{subsec:semicomp}.} $u,v\in C(\bar\Om)$ are semi-convex and semi-concave in $\Om$, $\LL u\le0\le\LL v$ in $\Om$ in the viscosity sense,
\begin{equation}\label{eq:bdtau}
u-v\le-\tau\ \text{on}\ \del\Om\ \text{for some}\ \tau>0,\quad\ M_0:=\max_{\bar\Om}(u-v)\geq\tau,\quad\ E:=\argmax_\Om(u-v)\neq\emp .
\end{equation}
Let $\om_u,\om_v$ be non-decreasing moduli of continuity of $u,v$ on $\bar\Om$ with respect to $d$ (they exist because $\bar\Om$ is compact and $|x-y|\le c\,d(x,y)$), and fix $\rho_0>0$ with
\begin{equation}\label{eq:rho0}
\om_u(\rho_0)+\om_v(\rho_0)\le\tau/4,\qquad\text{and set}\qquad \delta_1:=\rho_0/4 .
\end{equation}
Finally let $L_u,L_v$ be the Lipschitz constants of $u$ and $v$ on $K_{\rho_0/4}$ with respect to $d$; they are finite by Section~\ref{subsec:semi}(3) and depend on $\|u\|_\infty,\|v\|_\infty$, the semi-convexity constants, $\rho_0$, $\Om$ and $\G$.
\begin{Lem}[localization]\label{lem:bdry}
Under \eqref{eq:bdtau}--\eqref{eq:rho0}, let $0<\delta\le\delta_1$ and $\|h\|<\delta$. Then:
\begin{enumerate}[(a)]
\item $E\subset K_{\rho_0}$, i.e.\ $\dist(E,\del\Om)\ge\rho_0$; in particular $M_\delta(0)=M_0$ and $\A_\delta(0)=E$;
\item $u_h-v\le-\tau/2$ on $\Om_\delta\mns K_{\rho_0/2}$;
\item $M_\delta(h)\ge M_0-\tau/4>0$, the supremum in \eqref{eq:MAl} is attained, and $\emp\ne\A_\delta(h)\subset K_{\rho_0/2}$; moreover $\A_\delta(h)*h\subset K_{\rho_0/4}$.
\end{enumerate}
\end{Lem}
\begin{proof}
If $\dist(y,\del\Om)<\rho_0$, pick $y'\in\del\Om$ with $d(y,y')<\rho_0$; then $(u-v)(y)\le(u-v)(y')+\om_u(\rho_0)+\om_v(\rho_0)\le-3\tau/4<0<M_0$, which proves (a), since $\delta<\rho_0$. For (b), let $x\in\Om_\delta$ with $\dist(x,\del\Om)<\rho_0/2$; then $\dist(x*h,\del\Om)<\rho_0/2+\delta<\rho_0$, so $(u-v)(x*h)\le-3\tau/4$ as above, while $|v(x*h)-v(x)|\le\om_v(\|h\|)\le\tau/4$ by \eqref{eq:righttrans} and \eqref{eq:rho0}; hence $(u_h-v)(x)=(u-v)(x*h)+v(x*h)-v(x)\le-\tau/2$. For (c), take $y\in E$; by (a), $\dist(y*\inv h,\del\Om)\ge\rho_0-\delta>\delta$, so $y*\inv h\in\Om_\delta$, and $(u_h-v)(y*\inv h)=u(y)-v(y*\inv h)\ge M_0-\om_v(\|h\|)\ge M_0-\tau/4\geq 3\tau/4$. Since $M_0-\tau/4>-\tau/2$, by (b) the supremum in \eqref{eq:MAl} is a supremum over the compact set $K_{\rho_0/2}$, hence attained there; the last assertion follows from \eqref{eq:righttrans}.
\end{proof}
Thus, under \eqref{eq:bdtau}, all maximal sets $\A_\delta(h)$, $\delta\le\delta_1$, live in the fixed compact set $K_{\rho_0/2}\subset\Om$, on a neighbourhood of which $u$ and $v$ are Lipschitz. No Lipschitz bound of $u,v$ up to $\del\Om$ is needed anywhere: the boundary is handled by uniform continuity, and this is what makes the constants below independent of the regularization parameter in Section~\ref{subsec:compvisc}.
\begin{Lem}[Lipschitz dependence on the parameter]\label{lem:Mlip}
Under \eqref{eq:bdtau}--\eqref{eq:rho0}, for $0<\delta\le\delta_1$ and $h,h'\in B_\delta(0)$,
\begin{equation}\label{eq:lip}
|M_{\delta}(h)- M_{\delta}(h')|\leq L_u\, d(h,h').
\end{equation}
\end{Lem}
\begin{proof}
Let $x\in \A_{\delta}(h)$ and $x'\in \A_{\delta}(h')$. By maximality at $x'$ and Lemma~\ref{lem:bdry}(c),
\begin{align*}M_{\delta}(h)- M_{\delta}(h')  &= u(x*h)-v(x)-u(x'*h')+v(x')\\
&\le  u(x*h)-u(x*h')\le  L_u\,d(x*h,x*h')=L_u\,d(h,h'),\end{align*}
since $x*h,x*h'\in K_{\rho_0/4}$ and $d(x*h,x*h')=\|\inv{h'}*h\|=d(h,h')$. Exchanging the roles of $h$ and $h'$ gives \eqref{eq:lip}.
\end{proof}
The horizontal differential of $M_\delta$, where it exists, is determined by $\X u$ at the maximal points.
\begin{Lem}[envelope identity]\label{lem:envelope}
Under \eqref{eq:bdtau}--\eqref{eq:rho0}, let $0<\delta\le\delta_1$ and let $h\in B_\delta(0)$ be a point at which $M_\delta$ is Pansu differentiable, with horizontal differential $\X_hM_\delta(h)\in\R^m$. Then
$$\X u(x*h)=\X_hM_\delta(h)\qquad\text{for every}\ x\in\A_\delta(h);$$
in particular all maximal points share the same horizontal gradient. By Lemma~\ref{lem:Mlip} and Pansu's theorem this holds for almost every $h\in B_\delta(0)$.
\end{Lem}
\begin{proof}
Fix $x\in\A_\delta(h)$ and $\zeta\in\R^m$, and let $g_t=\exp\big(t\sum_i\zeta_iX_i\big)$, so that $\pi_1(g_t)=t\zeta$ and $\|g_t\|=|t||\zeta|$. For $|t|$ small, $h*g_t\in B_\delta(0)$ and, since $x\in\Om_\delta$, $M_\delta(h*g_t)\ge u(x*h*g_t)-v(x)=M_\delta(h)+t\inp{\X u(x*h)}{\zeta}+o(t)$ by \eqref{eq:taylor2}, $u$ being differentiable at $x*h$. Pansu differentiability of $M_\delta$ at $h$ gives $M_\delta(h*g_t)=M_\delta(h)+t\inp{\X_hM_\delta(h)}{\zeta}+o(t)$. Comparing for $t>0$ and $t<0$ yields $\inp{\X u(x*h)-\X_hM_\delta(h)}{\zeta}=0$ for every $\zeta$.
\end{proof}

We now come to the vanishing-gradient mechanism of \cite{Barles-Busca}, the only place where the strong maximum principle enters.
\begin{Lem}[non-degeneracy at every scale]\label{lem:case1}
Under \eqref{eq:bdtau}--\eqref{eq:rho0}, let $0<\delta\le\delta_1$ and $0<\eta\le\delta$. Then there exists $h\in B_\eta(0)$ such that
$$\X u(x*h)\ne0\qquad\text{for all}\ x\in\A_\delta(h).$$
\end{Lem}
\begin{proof}
Assume by contradiction that for every $h\in B_\eta(0)$ there is $x_h\in\A_\delta(h)$ with $\X u(x_h*h)=0$. By Lemma~\ref{lem:Mlip}, $M_\delta$ is $d$-Lipschitz on $B_\delta(0)$, so by Pansu's theorem it is Pansu differentiable at almost every $h\in B_\eta(0)$, and at every such $h$ Lemma~\ref{lem:envelope} gives $\X_hM_\delta(h)=\X u(x_h*h)=0$. Since $B_\eta(0)$ is connected, $M_\delta$ is constant on $B_\eta(0)$ \cite[Proposition 4.8]{DMV13}, and by Lemma~\ref{lem:bdry}(a), $M_\delta(h)=M_\delta(0)=M_0$ for all $h\in B_\eta(0)$.

Now let $y\in E$ and $\|g\|<\eta$. Then $y\in\Om_\delta$ by Lemma~\ref{lem:bdry}(a), and
$$u(y*g)-v(y)=u_g(y)-v(y)\le M_\delta(g)=M_0=u(y)-v(y),$$
that is, $u\le u(y)$ on $B_\eta(y)=y*B_\eta(0)$, a ball contained in $\Om$ because $\eta\le\delta<\rho_0$. By Corollary~\ref{cor:stmx} applied to the sub-solution $u-u(y)$ on $B_\eta(y)$, $u\equiv u(y)$ on $B_\eta(y)$. Then, for $y'\in B_\eta(y)$, $(u-v)(y')\le M_0=(u-v)(y)$ and $u(y')=u(y)$ give $v(y')\ge v(y)$: the super-solution $v-v(y)$ attains a non-positive minimum at $y$ on $B_\eta(y)$, and by Corollary~\ref{cor:stmx} again $v\equiv v(y)$ on $B_\eta(y)$. Therefore $B_\eta(y)\subset E$ for every $y\in E$: $E$ is open. It is also relatively closed in $\Om$ and non-empty, so it contains a connected component $C$ of $\Om$ (we do not assume $\Om$ connected, because reduced domains will be used in Section~\ref{subsec:compvisc}); thus $u-v\equiv M_0>0$ on $C$, while $\del C\subset\del\Om$ where $u-v\le-\tau<0$, which contradicts the continuity of $u-v$ on $\bar\Om$.
\end{proof}
Lemma~\ref{lem:case1} is the Carnot group version of the alternative in the second step of the Euclidean argument (Section~\ref{subsec:BB}). In $\R^n$ it is all that is needed, since $u_h$ is then a sub-solution (Theorem~\ref{thm:uncond}(iii) below). In general $u_h$ is not a sub-solution (Proposition~\ref{prop:central}), and the lemma provides $\X u\neq0$ on $\A_\delta(h)*h$ rather than on $E$; see Section~\ref{sec:obstruction}.

\subsection{Comparison for viscosity solutions}\label{subsec:compvisc}
Viscosity sub/super-solutions can be approximated by sup/inf-convolutions, which are semi-convex/semi-concave functions \cite{Ar-Cr-Ju, Cr-Ish-Lions}. In Carnot groups this was done by Wang \cite{CYW07} with a kernel built from the group structure. Set $q:=2\,r!$, so that $\|x\|^q$ is a polynomial, and for $w\in C(\bar \Om)$ and $\eps>0$ let
\begin{equation}\label{eq:supinfconv}
w^\eps(x)= \sup_{y\in \bar \Om} \Big\{w(y) -\frac{1}{2\eps}\, \|y*\inv x\|^{q}\Big\},\qquad
w_\eps(x)= \inf_{y\in \bar \Om} \Big\{w(y) +\frac{1}{2\eps}\, \|y*\inv x\|^{q}\Big\},\qquad x\in\bar\Om,
\end{equation}
where $\|y*\inv x\|=d(\inv x,\inv y)$ is the right-invariant distance. Substituting $y=z*x$,
\begin{equation}\label{eq:leftsup}
w^\eps(x)= \sup\Big\{w(z*x) -\frac{\|z\|^{q}}{2\eps}:\ z*x\in\bar\Om\Big\},
\end{equation}
so that $w^\eps$ is a supremum of \emph{left} translates of $w$; this is what makes \eqref{eq:supinfconv} compatible with the equation. For $R:=2\max\{\|u\|_{L^\infty},\|v\|_{L^\infty}\}$ let
\begin{equation}\label{eq:redom}
\Om[\eps]:=\Big\{x\in \Om: \inf_{y\in \R^n\mns \Om}\|y*\inv x\|^{q} > (1+4R)\eps\Big\},
\end{equation}
a bounded open set with $\Om[\eps]\uparrow\Om$ as $\eps\to0^+$, and $\sup_{z\in\del\Om[\eps]}\dist(z,\del\Om)\to0$ as $\eps\to0^+$ (by continuity of the group operations on bounded sets).
\begin{Lem}\label{lem:wang}
Let $u,v\in C(\bar\Om)$ and $\eps>0$. Then:
\begin{enumerate}
\item $u^\eps$ is semi-convex and $v_\eps$ is semi-concave in $\Om$;
\item $u^\eps\downarrow u$ and $v_\eps\uparrow v$ uniformly on $\bar\Om$ as $\eps\to0^+$;
\item if $\LL u\leq 0\leq \LL v$ in $\Om$ in the viscosity sense, then $\LL u^\eps\leq 0\leq \LL v_\eps$ in $\Om[\eps]$ in the viscosity sense;
\item if $\om_u$, $\om_v$ are moduli of continuity of $u,v$ on $\bar\Om$ with respect to $d$, then $|u^\eps(x)-u^\eps(x')|\le\om_u(d(x,x'))$ and $|v_\eps(x)-v_\eps(x')|\le\om_v(d(x,x'))$ for all $x,x'\in\Om[\eps]$; moreover, at every point $x\in\Om[\eps]$ the supremum in \eqref{eq:leftsup} is attained at some $z_x$ with $\|z_x\|^q\le 4\eps\|u\|_{L^\infty}$, and $z_x*x'\in\Om$ for every $x'\in\Om[\eps]$; the same holds for the infimum $v_\eps(x)=\inf\{v(z*x)+\|z\|^q/2\eps:\ z*x\in\bar\Om\}$.
\end{enumerate}
\end{Lem}
\begin{proof}
(1)--(3) are \cite[Proposition 3.3]{CYW07}: $\|y*\inv x\|^q$ is a polynomial in $x$, so $u^\eps$ is a supremum of smooth functions with Hessians bounded below by $-C/\eps$ on bounded sets, and the proof of (3) transports a test function touching $u^\eps$ at $x$ to the test function $\varphi\circ\inv{L_{z_x}}+\text{const}$ touching $u$ at $z_x*x$, which has the same horizontal gradient by left invariance; the class \eqref{eq:defA} is preserved. For (4), let $x\in\Om[\eps]$; the supremum in \eqref{eq:supinfconv} is attained at some $y_x\in\bar\Om$, and comparing with $y=x$ gives $\|y_x*\inv x\|^q\le2\eps\,(u(y_x)-u(x))\le4\eps\|u\|_\infty<(1+4R)\eps$; with $z_x:=y_x*\inv x$, for $x'\in\Om[\eps]$ the point $z_x*x'$ satisfies $\|(z_x*x')*\inv{x'}\|^q=\|z_x\|^q<(1+4R)\eps$, hence lies in $\Om$ by \eqref{eq:redom}. Therefore $u^\eps(x')\ge u(z_x*x')-\|z_x\|^q/2\eps$ and
$$u^\eps(x)-u^\eps(x')\le u(z_x*x)-u(z_x*x')\le\om_u\big(d(z_x*x,z_x*x')\big)=\om_u(d(x,x')),$$
by left invariance of $d$. Exchanging $x$ and $x'$ proves the claim for $u^\eps$; $v_\eps$ is treated in the same way.
\end{proof}
Property (4) is the point of using the kernel \eqref{eq:supinfconv}: the regularized functions inherit the moduli of continuity of $u$ and $v$ uniformly in $\eps$, whereas their Lipschitz constants blow up as $\eps\to0^+$ unless $u,v$ are Lipschitz.

\begin{proof}[Proof of Theorem \ref{thm:comp}]
Suppose $M_0:=\max_{\bar\Om}(u-v)>0$ and let $\tau:=M_0/2$, $\tilde v:=v+\tau$, so that $u-\tilde v\le-\tau$ on $\del\Om$ and $\max_{\bar\Om}(u-\tilde v)=\tau>0$; note that $\tilde v_\eps=v_\eps+\tau$, so that $E_\eps$ in \eqref{eq:Eeps} is also the set of points of $\Om[\eps]$ at which $u^\eps-\tilde v_\eps$ attains its maximum over $\overline{\Om[\eps]}$. Let $\om_u,\om_v$ be moduli of continuity of $u,v$ on $\bar\Om$ with respect to $d$, and fix $\rho_0>0$ with $\om_u(\rho_0)+\om_v(\rho_0)\le\tau/8$. The set $E=\argmax_\Om(u-v)$ is compact and $\dist(E,\del\Om)\ge\rho_0$: if $\dist(y,\del\Om)<\rho_0$ and $y'\in\del\Om$ satisfies $d(y,y')<\rho_0$, then $(u-\tilde v)(y)\le(u-\tilde v)(y')+\tau/8\le-7\tau/8<\tau$. By Lemma~\ref{lem:wang}(2) and the properties of $\Om[\eps]$ noted after \eqref{eq:redom}, there is $\eps_1>0$ such that for $\eps<\eps_1$: $\|u^\eps-u\|_\infty+\|v_\eps-v\|_\infty\le\tau/8$, $\sup_{z\in\del\Om[\eps]}\dist(z,\del\Om)<\rho_0$, and $\Om[\eps]\supset E$. Fix $\eps=\eps_k<\eps_1$ and consider the pair $(u^\eps,\tilde v_\eps)$ on $\Om[\eps]$. It is semi-convex/semi-concave and continuous on $\overline{\Om[\eps]}$, $\LL u^\eps\le0\le\LL\tilde v_\eps$ in $\Om[\eps]$ (Lemma~\ref{lem:wang}(1),(3)), and:
\begin{itemize}
\item for $z\in\del\Om[\eps]$ there is $z'\in\del\Om$ with $d(z,z')<\rho_0$, so $(u-\tilde v)(z)\le-\tau+\tau/8$ and $(u^\eps-\tilde v_\eps)(z)\le(u-\tilde v)(z)+\tau/4\le-\tau/2$;
\item $\max_{\overline{\Om[\eps]}}(u^\eps-\tilde v_\eps)\ge\max_E(u-\tilde v)=\tau>0$, since $u^\eps\ge u$ and $\tilde v_\eps\le\tilde v$.
\end{itemize}
Hence the maximum of $u^\eps-\tilde v_\eps$ over $\overline{\Om[\eps]}$ is attained only in $\Om[\eps]$. Thus $E_\eps=\argmax_{\Om[\eps]}(u^\eps-\tilde v_\eps)$ is a non-empty compact subset of $\Om[\eps]$ and $u^\eps\le\tilde v_\eps$ on $\del\Om[\eps]$. In other words, the pair $(u^\eps,\tilde v_\eps)$ satisfies the standing assumptions \eqref{eq:bdtau} on $\Om[\eps]$, with $\tau/2$ in place of $\tau$. By hypothesis $\X u^\eps\ne0$ on $E_\eps$, so Proposition~\ref{prop:comp0}, applied to $(u^\eps,\tilde v_\eps)$ on $\Om[\eps]$, gives $u^\eps\le\tilde v_\eps$ in $\Om[\eps]$, contradicting the second item.
\end{proof}

\begin{proof}[Proof of Theorem \ref{thm:uncond}]
(i) Let $N$ be an open neighbourhood of $E=\argmax_\Om(u-v)$ on which $v$ is $C^1$ with $\X v\ne0$ (the case of $u$ is symmetric, with $u^\eps$ in place of $v_\eps$). If $M_0\le0$ there is nothing to prove, so suppose $M_0>0$ and use the notation of the previous proof. For $\eps<\eps_1$ the set $E_\eps$ is a non-empty compact subset of $\Om[\eps]$, and $\sup_{E_\eps}\dist(\cdot,E)\to0$ as $\eps\to0^+$: indeed $u^\eps-\tilde v_\eps\to u-\tilde v$ uniformly on $\bar\Om$ and $(u^\eps-\tilde v_\eps)\ge\tau$ on $E_\eps$, so every limit point of points of $E_\eps$ is a maximum point of $u-\tilde v$ over $\bar\Om$, hence belongs to $E$. Let $x\in E_\eps$; then $v_\eps$ is differentiable at $x$ (Section~\ref{subsec:semi}(1)). By Lemma~\ref{lem:wang}(4) applied to $v_\eps$, there is $z_x$ with $\|z_x\|^q\le4\eps\|v\|_\infty$ and $z_x*x'\in\Om$ for all $x'\in\Om[\eps]$ such that
$$v_\eps(x')\le v(z_x*x')+\frac{\|z_x\|^q}{2\eps}\quad\text{for}\ x'\in\Om[\eps],\qquad\text{with equality at}\ x'=x .$$
By \eqref{eq:trinangle}, $d(z_x*x,x)=\|\inv x*z_x*x\|\le c_\G\|z_x\|^{1/r}\to0$ uniformly in $x$ as $\eps\to0^+$, so for $\eps$ small the points $z_x*x$, $x\in E_\eps$, lie in $N$, and so do the points $z_x*x'$ for $x'$ near $x$. The right-hand side is then a $C^1$ function of $x'$ near $x$ touching the differentiable function $v_\eps$ from above at $x$, so their gradients coincide at $x$, and by left invariance $\X v_\eps(x)=\X(v\circ L_{z_x})(x)=(\X v)(z_x*x)\ne0$. Since $\X u^\eps=\X v_\eps$ on $E_\eps$, the hypothesis of Theorem~\ref{thm:comp} holds for every small $\eps$, and Theorem~\ref{thm:comp} gives $u\le v$, against $M_0>0$.

(ii) If $v\equiv c$ and $M_0>0$, then $E=\argmax_\Om(u-v)=\argmax_\Om u$ is a non-empty compact subset of $\Om$ and $u$ attains an interior maximum; adding a constant we may take it non-negative, and Corollary~\ref{cor:stmx} gives that $u$ is constant on the component of $\Om$ containing a point of $E$, so $u-v\equiv M_0>0$ there, contradicting $u\le v$ on $\del\Om$. The case $u$ constant uses the strong minimum principle for $v$.

(iii) Suppose first $r=1$, so that $\G=(\R^n,+)$, $L_h=R_h$, and $d$ is a norm distance. By the proof of Theorem~\ref{thm:comp}, it suffices to show that under the standing assumptions \eqref{eq:bdtau}--\eqref{eq:rho0} a semi-convex/semi-concave pair cannot exist. Given such a pair, apply Lemma~\ref{lem:case1} with $\eta=\delta\le\delta_1$: there is $h\in B_\delta(0)$ with $\X u(\cdot+h)\ne0$ on $\A_\delta(h)=\argmax_{\Om_\delta}(u_h-v)$. Since $u_h=u^h$ is a sub-solution in $\Om_\delta$ by Lemma~\ref{lem:lefttrans}, $u_h-v\le-\tau/2$ on $\del\Om_\delta$ by Lemma~\ref{lem:bdry}(b) and continuity, and $M_\delta(h)>0$ is attained in $\Om_\delta$, Proposition~\ref{prop:comp0} applied to $(u_h,v)$ on $\Om_\delta$ gives $u_h\le v$ on $\Om_\delta$, a contradiction. (This is the argument of \cite{Barles-Busca}: in the Euclidean case the translation that carries the Taylor expansion is also a symmetry of the equation.)

Now let $r\ge1$ and $u=F\circ\pi_1$, $v=G\circ\pi_1$ with $F,G$ continuous on $\pi_1(\bar\Om)$. The set $\Om':=\pi_1(\Om)\subset\R^m$ is open and bounded, and $\del\Om'\subset\pi_1(\del\Om)$: if $\zeta\in\del\Om'$, $\zeta=\lim\pi_1(x_k)$ with $x_k\in\Om$, and a limit point $\bar x\in\bar\Om$ of $(x_k)$ has $\pi_1(\bar x)=\zeta$ and cannot lie in the open set $\Om$ (otherwise $\zeta\in\Om'$). Hence $F\le G$ on $\del\Om'$. In exponential coordinates $X_i=\del_{x_{i,1}}+\sum_{j\ge2}\sum_kb^i_{k,j}\del_{x_{k,j}}$, so for $\psi\in C^2(\Om')$ and $\varphi:=\psi\circ\pi_1$ one has $X_i\varphi=(\del_i\psi)\circ\pi_1$ and $X_jX_i\varphi=(\del_j\del_i\psi)\circ\pi_1$; thus $\LL\varphi=(\LL^{\rm eucl}\psi)\circ\pi_1$, where $\LL^{\rm eucl}\psi=-\Tr(A(\gr\psi)D^2\psi)+H(\gr\psi)$ is the Euclidean operator with the same $A,H$. If $\psi$ touches $F$ from above at $\zeta\in\Om'$ with $\gr\psi(\zeta)\ne0$, then $\varphi$ touches $u$ from above at every $x\in\Om$ with $\pi_1(x)=\zeta$, with $\X\varphi(x)=\gr\psi(\zeta)\ne0$, so $\LL^{\rm eucl}\psi(\zeta)=\LL\varphi(x)\le0$: $F$ is a viscosity sub-solution of $\LL^{\rm eucl}F=0$ in $\Om'$, and likewise $G$ is a super-solution. By the case $r=1$ applied in $\R^m$, $F\le G$ in $\Om'$, i.e.\ $u\le v$ in $\Om$.
\end{proof}

\section{The obstruction and an open problem}\label{sec:obstruction}
In this section we make precise the statement, used informally in the introduction, that the non-commutativity of $\G$ separates the two roles of the Euclidean translation, we quantify what is missing to remove the hypothesis of Theorem~\ref{thm:comp}, and we formulate the open problem.

\subsection{Right translations and the equation}
For $g\in\G$ and a smooth $\varphi$, the chain rule gives
\begin{equation}\label{eq:pushforward}
X_i(\varphi\circ R_g)=(\tilde X_i\varphi)\circ R_g,\qquad \tilde X_i:=(R_g)_*X_i=\operatorname{Ad}_{\inv g}X_i,
\end{equation}
where $(R_g)_*X_i$ is the push-forward of $X_i$ by the diffeomorphism $R_g$; it is a left-invariant vector field because $R_g$ commutes with left translations, and it equals $\operatorname{Ad}_{\inv g}X_i=X_i+Z_i$ with $Z_i\in\g_2\oplus\ldots\oplus\g_r$ (write $g=\exp(Y)$ and expand $\operatorname{Ad}_{\inv g}=e^{-\operatorname{ad}Y}$: all terms beyond the first lie in the higher layers). Consequently
\begin{equation}\label{eq:Ltilde}
\LL(u\circ R_g)(x)=\tilde\LL u(x*g),\qquad \tilde\LL u:=-\Tr\big(A(\tilde\X u)\tilde\X\tilde\X^\star u\big)+H(\tilde\X u),
\end{equation}
where $\tilde\LL$ is the operator \eqref{eq:maineq} built with the fields $\tilde X_1,\ldots,\tilde X_m$. Let $\Z(\G)$ denote the center of $\G$; it contains $\exp(\g_r)$.
\begin{Prop}\label{prop:central}
Let $g\in\G$.
\begin{enumerate}[(i)]
\item $\tilde X_i=X_i$ for all $i=1,\ldots,m$, i.e.\ $X_i(\varphi\circ R_g)=(X_i\varphi)\circ R_g$ for all smooth $\varphi$, if and only if $g\in\Z(\G)$.
\item If $g\in\Z(\G)$, then $R_g=L_g$, so $u\circ R_g$ is a viscosity sub-solution (super-solution) on $\inv g*\Om$ whenever $u$ is one on $\Om$ (Lemma~\ref{lem:lefttrans}); and for a one-parameter subgroup $g_t=\exp(tY)$ with $Y\in\g$ central (for instance $Y\in\g_r$), $\frac{d}{dt}u(x*g_t)=Yu(x)$ for $u\in C^1$, a derivative of $u$ along a non-horizontal direction.
\item If $g\notin\Z(\G)$, right translation by $g$ need not preserve sub-solutions. For instance, in the Heisenberg group $\hh^1$ with $X_1=\del_1-\frac12x_2\del_3$, $X_2=\del_2+\frac12x_1\del_3$ and the $\infty$-Laplacian, the function $u(y)=y_1y_3$ satisfies $\Delta_{\X,\infty}u(2,-1,0)=5>0$, hence is a classical, strict sub-solution of $-\Delta_{\X,\infty}u=0$ in a neighbourhood $U$ of $(2,-1,0)$, while for $g=(0,-1,0)$ one has $\tilde X_1=X_1-\del_3$, $\tilde X_2=X_2$, $\tilde\X u(2,-1,0)=(-1,2)\ne0$ and $\tilde\Delta_{\X,\infty}u(2,-1,0)=-5<0$: by \eqref{eq:Ltilde}, $u\circ R_g$ is not a sub-solution on $\inv{R_g}(U)$. By the dilations and the rotations of $\hh^1$, which preserve the $\infty$-Laplacian, the same happens for every non-central $g$.
\end{enumerate}
\end{Prop}
\begin{proof}
(i) $\tilde X_i=X_i$ for all $i$ iff $\operatorname{Ad}_{\inv g}$ fixes $\g_1$, iff $g$ commutes with $\exp(tX_i)$ for all $t\in\R$ and $i$, iff $g$ commutes with the subgroup generated by the horizontal one-parameter subgroups, which is $\G$ by \eqref{eq:horcond} (Chow's theorem). (ii) is immediate from $x*g_t=g_t*x$ and the definition of $Yu$. For (iii) one computes: $X_1u=y_3-\frac12y_1y_2$, $X_2u=\frac12y_1^2$, $\XXs u=\begin{pmatrix}-y_2&\frac12y_1\\[2pt] \frac12y_1&0\end{pmatrix}$, so $\Delta_{\X,\infty}u=-y_2\big(y_3-\frac12y_1y_2\big)^2+\frac12y_1^3\big(y_3-\frac12y_1y_2\big)$, which equals $5$ at $(2,-1,0)$; and for $g=(0,b,0)$ we have $\tilde X_1=X_1+b\,\del_3$ and $\tilde X_2=X_2$, so that, with $F:=\tilde X_1u=y_3-\frac12y_1y_2+by_1$, $\tilde X_1\tilde X_1u=2b-y_2$, $\tilde X_2\tilde X_2u=0$, $\tfrac12(\tilde X_1\tilde X_2+\tilde X_2\tilde X_1)u=\frac12y_1$, and $\tilde\Delta_{\X,\infty}u=(2b-y_2)F^2+\frac12y_1^3F$, which at $(2,-1,0)$ and $b=-1$ equals $-1-4=-5$.
\end{proof}

Thus the two properties of the Euclidean translation used in the Barles--Busca argument are shared, on a non-abelian Carnot group, only by the central translations, and along those the parameter derivative of $t\mapsto u(x*g_t)$ is a non-horizontal derivative of $u$, which the equation does not control. The vanishing of $Yu$ at the maximal points of $u-v$ yields a maximum of $u$ along the curve $t\mapsto x*g_t$, not on an open set, so the strong maximum principle does not apply. This is not an artifact of the argument: on any Carnot group the function $u=x_1^2$ (the square of a horizontal coordinate) is a smooth, non-constant sub-solution of the $\infty$-Laplace equation and of \eqref{eq:diag} ($\X u=(2x_1,0,\ldots,0)$, $\XXs u=\operatorname{diag}(2,0,\ldots,0)$, $\LL u=-8x_1^2\le0$), it is invariant under the center, and its horizontal gradient vanishes on the hypersurface $\{x_1=0\}$.

\subsection{The open problem}\label{subsec:openpb}
By the proof of Theorem~\ref{thm:comp}, the comparison principle for all continuous sub- and super-solutions is equivalent to the statement that no semi-convex/semi-concave pair satisfies the standing assumptions \eqref{eq:bdtau}--\eqref{eq:rho0}.
\begin{Prob}\label{prob:main}
Let $\G$ be a non-abelian Carnot group and let $A$, $H$ satisfy \eqref{eq:ell} and \eqref{eq:scale}. Does the comparison principle hold for all $u,v\in C(\bar\Om)$ with $\LL u\le0\le\LL v$ in $\Om$ and $u\le v$ on $\del\Om$, that is, can the hypothesis of Theorem~\ref{thm:comp} be removed? Equivalently: can the standing assumptions \eqref{eq:bdtau}--\eqref{eq:rho0} be satisfied by a semi-convex/semi-concave pair?
\end{Prob}
For the $\infty$-Laplacian the answer is positive by the approximation method, and for the $p$-Laplacian through weak solutions (Section~\ref{subsec:examples}); for the other equations of Table~\ref{tab:examples} that satisfy (E) and (S) the question is, to our knowledge, open.

The Euclidean proof settles the problem by Lemma~\ref{lem:case1} and Proposition~\ref{prop:comp0} applied to the translate $(u_h,v)$. On a non-abelian group $u_h$ is not a sub-solution unless $h$ is central (Proposition~\ref{prop:central}), Proposition~\ref{prop:comp0} can only be applied to $(u,v)$ itself, which requires $\X u\ne0$ on $E=\A_\delta(0)$, and Lemma~\ref{lem:case1} provides $\X u\ne0$ on $\A_\delta(h)*h$ for some $h$ in every ball around $0$. The following proposition shows what would be needed to pass from one to the other. For $h\in B_\delta(0)$ let
\begin{equation}\label{eq:thetarho}
\theta(h):=\min_{x\in\A_\delta(h)}|\X u(x*h)|,\qquad \rho_h(\eps):=\sup\big\{\dist(z,\A_\delta(h)):\ z\in\A_\delta^\eps(h)\big\},
\end{equation}
where $\A_\delta^\eps(h):=\{z\in\Om_\delta:\ (u_h-v)(z)\ge M_\delta(h)-\eps\}$ is the set of $\eps$-approximate maximizers of $u_h-v$, so that $\rho_h(\eps)$ is its excess over $\A_\delta(h)$, and let
\begin{equation}\label{eq:modhl}
\om_h(t):=\sup\Big\{|\X u(y)-\X u(\hat y)|:\ \hat y\in\A_\delta(h)*h,\ y\in K_{\rho_0/4}\cap D_u,\ d(y,\hat y)\le t\Big\},
\end{equation}
where $D_u$ is the set of differentiability points of $u$, be the modulus of continuity of $\X u$ relative to the compact set $\A_\delta(h)*h\subset D_u$. (For a general semi-convex function no modulus of continuity of the gradient on its set of differentiability points exists: $u=|x|$ on $\R$ has gradients $\pm1$ at the points $\pm s$. What exists, by partial continuity, is a modulus relative to a compact set of differentiability points.) By Lemma~\ref{lem:envelope}, $\theta(h)=|\X_hM_\delta(h)|$ at the points of Pansu differentiability of $M_\delta$.
\begin{Prop}[transfer of non-degeneracy]\label{prop:transfer}
Under \eqref{eq:bdtau}--\eqref{eq:rho0}, let $0<\delta\le\delta_1$ and $h\in B_\delta(0)$. Then $\rho_h(\eps)\to0$ as $\eps\to0^+$, $\om_h(t)\to0$ as $t\to0^+$, and
\begin{equation}\label{eq:transfer}
|\X u(y)|\ \ge\ \theta(h)-\om_h\big(\rho_h(2L_u\|h\|)+\|h\|\big)\qquad\text{for every}\ y\in E .
\end{equation}
In particular, if there are $h_k\to0$ with $\theta(h_k)\ge c>0$ and $\om_{h_k}\big(\rho_{h_k}(2L_u\|h_k\|)+\|h_k\|\big)\to0$, then $\X u\ne0$ on $E$, and Proposition~\ref{prop:comp0} shows that the pair $(u,v)$ cannot exist.
\end{Prop}
\begin{proof}
If $\rho_h(\eps_k)\ge\rho>0$ along a sequence $\eps_k\to0$, pick $z_k\in\Om_\delta$ with $(u_h-v)(z_k)\ge M_\delta(h)-\eps_k$ and $\dist(z_k,\A_\delta(h))\ge\rho/2$. For $\eps_k\le\tau/4$ we have $(u_h-v)(z_k)\ge M_0-\tau/2>-\tau/2$, so $z_k\in K_{\rho_0/2}$ by Lemma~\ref{lem:bdry}(b),(c), and a subsequence converges to some $z\in K_{\rho_0/2}$ with $(u_h-v)(z)=M_\delta(h)$ and $\dist(z,\A_\delta(h))\ge\rho/2$, a contradiction. If $\om_h(t_k)\ge\eps_0>0$ along $t_k\to0$, there are $\hat y_k\in\A_\delta(h)*h$ and $y_k\in K_{\rho_0/4}\cap D_u$ with $d(y_k,\hat y_k)\le t_k$ and $|\X u(y_k)-\X u(\hat y_k)|\ge\eps_0$; along a subsequence $\hat y_k\to\hat y\in\A_\delta(h)*h$, hence $y_k\to\hat y$, and partial continuity of the gradient at the differentiability point $\hat y$ gives $\gr u(y_k)\to\gr u(\hat y)$ and $\gr u(\hat y_k)\to\gr u(\hat y)$, a contradiction since $\X u=\sigma^T\gr u$ with $\sigma$ continuous.

Let $y\in E$. By Lemma~\ref{lem:bdry}(a), $y\in\Om_\delta$ and $y,y*h\in K_{\rho_0/4}$, so $(u_h-v)(y)=(u-v)(y)+u(y*h)-u(y)\ge M_0-L_u\|h\|$, while $M_\delta(h)\le M_\delta(0)+L_u\|h\|=M_0+L_u\|h\|$ by Lemma~\ref{lem:Mlip}. Hence $y$ is a $2L_u\|h\|$-approximate maximizer of $u_h-v$, and by compactness there is $\hat x\in\A_\delta(h)$ with $d(y,\hat x)\le\rho_h(2L_u\|h\|)$. By \eqref{eq:righttrans}, $d(y,\hat x*h)\le d(y,\hat x)+d(\hat x,\hat x*h)\le\rho_h(2L_u\|h\|)+\|h\|$. Since $y$ is a differentiability point of $u$ in $K_{\rho_0/4}$ and $\hat x*h\in\A_\delta(h)*h$, the definition of $\om_h$ gives $|\X u(y)-\X u(\hat x*h)|\le\om_h\big(\rho_h(2L_u\|h\|)+\|h\|\big)$, while $|\X u(\hat x*h)|\ge\theta(h)$.
\end{proof}
A positive answer to Problem~\ref{prob:main} along the lines of Proposition~\ref{prop:transfer} requires two uniformity properties as $h\to0$ that we are not able to prove: a lower bound $\theta(h_k)\ge c$ along a sequence $h_k\to0$, that is, by Lemma~\ref{lem:envelope}, a sequence of points of Pansu differentiability at which the horizontal differential of the Lipschitz function $M_\delta$ does not degenerate (Lemma~\ref{lem:case1} gives $\theta(h)>0$ for some $h$ in every ball around $0$, but no lower bound), and the uniform control of the localization $\rho_{h_k}$ of the approximate maximizers and of the moduli $\om_{h_k}$. The strong maximum principle alone cannot give the first: if $\theta(h)\to0$ as $h\to0$ along the points of differentiability, one only obtains $M_\delta(h)=M_0+o(\|h\|)$ and, at the maximal points $y$ of $u-v$, the almost-maximality $\sup_{B_\eta(y)}u\le u(y)+o(\eta)$ for small $\eta$, and this is realized by the non-constant sub-solution $x_1^2$ at any point of $\{x_1=0\}$, where $\sup_{B_\eta}u-u(y)$ is of order $\eta^2$. Any proof has to use the function $M_\delta$ on a whole ball together with the boundary inequality, that is, it has to contain the same ingredients as the comparison principle itself. One may also allow a second parameter and consider $M_\delta(h,l):=\max_{\Om_\delta}\big(u(\cdot*h)-v(\cdot*l)\big)$: Lemma~\ref{lem:case1} and Proposition~\ref{prop:transfer} have analogues for the slices $h\mapsto M_\delta(h,l)$, with $E*\inv l$ in place of $E$, and it would suffice to find the uniform non-degeneracy near one $l$; we do not know whether this additional freedom helps. In the Euclidean case no transfer is needed, because Proposition~\ref{prop:comp0} applies to $(u_h,v)$ directly (Theorem~\ref{thm:uncond}(iii)), and for pairs invariant under $[\G,\G]$ the drift $\tilde X_i-X_i$ of \eqref{eq:pushforward} acts trivially, which is why the difficulty disappears there.
\section*{AI Statement}

While preparing the first version of this manuscript for submission, we asked an LLM (Claude) for a pre-submission review. The LLM identified a number of misprints, obscure points, and gaps. We were able to fix all but one, rendering our main result incomplete.  We are thankful to the LLM for detecting an error prior to publication. This version of the manuscript details the partial results that we are able to obtain at this time. It was completed with the help of the LLM. 

We have independently verified, validated, and rewritten all parts of the paper based on LLM suggestions and take full responsibility for the mathematical content of the paper. In addition, we have produce a formalization in Lean of the key results of this manuscript, \cite{CompND-Lean}. 
\section*{Acknowledgements}
J. Manfredi has been supported by the Simons collaboration grant 962828. S. Mukherjee has been supported by EPSRC New Investigator award [grant number EP/W001586/1]. \par
This work began when S. Mukherjee was supported by European Union's Horizon 2020 research and innovation program under the (GHAIA) Marie Sk\l{}odowska-Curie grant agreement No. 777822 during his visit to the University of Pittsburgh, and he thanks the Department of Mathematics there for their hospitality.


\bibliographystyle{plain}
\bibliography{MyBib}

\end{document}